\documentclass[10pt]{article}

\usepackage{graphicx} 
\usepackage{amsmath} 
\usepackage{amssymb} 
\usepackage{amsthm}
\usepackage{yfonts}
\usepackage{mathrsfs}
\usepackage{xypic}
\usepackage{bm}
\title{Injectivity radius for non-simply connected symmetric spaces via Cartan polyhedron } 
\author{Ling Yang\footnote{This work was partially supported by NNSFC.}}

\begin{document}
\def\e{\mathbf{e}}
\def\r{\mathbf{r}}
\def\om{\omega}
\def\Om{\Omega}
\def\td{\tilde}
\def\w{\wedge}
\def\c{\cdot}
\def\n{\nabla}
\def\y{\mathbf{y}}
\def\a{\mathbf{a}}
\def\b{\mathbf{b}}
\def\Y{\mathbf{Y}}
\def\p{\partial}
\def\f{\frac}
\def\si{\sigma}
\def\mc{\mathcal}
\def\C{\Bbb{C}}
\def\ra{\rightarrow}
\def\lan{\langle}
\def\ran{\rangle}
\def\i{\sqrt{-1}}
\def\la{\lambda}
\def\tr{\mbox{tr}}
\def\R{\Bbb{R}}
\def\Z{\Bbb{Z}}
\def\Q{\Bbb{Q}}
\def\ol{\overline}
\def\th{\theta}
\def\td{\tilde}
\def\e{\eta}
\def\ep{\epsilon}
\def\De{\Delta}
\def\a{\alpha}
\def\be{\beta}
\def\z{\zeta}
\def\La{\Lambda}
\def\la{\lambda}
\def\g{\gamma}
\def\de{\delta}
\def\Th{\Theta}
\def\p{\partial}
\def\k{\kappa}
\def\rank{\mbox{rank}}
\def\Id{\mathbf{Id}}
\def\d{\dot}
\def\dd{\ddot}
\def\sw{\textswab}
\def\Inn{\mbox{Inn}}
\def\wtd{\widetilde}
\def\Si{\Sigma}
\def\G{\Gamma}
\def\P{\Bbb{P}}
\def\H{\Bbb{H}}
\def\ze{\zeta}
\def\Exp{\mbox{Exp}}
\newtheorem{Pro}{Proposition}[section]
\newtheorem{Lem}{Lemma}[section]
\newtheorem{Thm}{Theorem}[section]
\newtheorem{Cor}{Corollary}[section]

\renewcommand{\theequation}{\arabic{section}.\arabic{equation}}

\maketitle

\begin{abstract}
We determine the cut locus of arbitrary non-simply connected, compact and irreducible Riemannian symmetric space
explicitly, and compute injectivity radius and diameter for
every type of them.
\end{abstract}

\section{Introduction}

Let $(M,g)$ be a Riemannian manifold, $p\in M$ and $\ze:[0,\infty)\ra M$ be a normal geodesic such that $\ze(0)=p$,
then the set of $t$ for with $d\big(\ze(t),\ze(0)\big)=t$ is either $[0,\infty)$ or $[0,t_0]$ for some $t_0>0$,
where $d(,)$ is the distance function on $M\times M$ induced by the metric $g$.
In the latter case, $\ze(t_0)$ is called the \textit{cut point of} $\ze$ with respect to $p$ and $t_0\dot{\ze}(0)$
is called a \textit{cut point in $T_p M$}. The union of all cut points in $T_p M$ is called the
\textit{cut locus of $p$ in $T_p M$} and denoted by $C(p)$. The \textit{injectivity radius}
of $M$ is the largest $r$ such that for all $p\in M$, $\exp_p$ is an embedding on the open ball of radius $r$
in $T_p M$, which is denoted by $i(M)$; the \textit{diameter} of $M$ is the least upper bound of the length of minimal
geodesics in $M$, which is denoted by $d(M)$.

$C(p),i(M),d(M)$ have a close relationship with other geometrical quantities, e.g., sectional curvature,
Ricci curvature, volume, fundamental group, conjugate locus, convexity radius and so on.
Cheeger, Klingenberg, Toponogov, Berger, Grove, Shiohama, Weinstein, Sugahara, Ichida and P\"{u}ttmann have made a contribution to these topics (see \cite{CE} Ch.5-6,
\cite{GS}, \cite{W}, \cite{Su}, \cite{Ic}, \cite{P}).

Generally, it is very difficult to determine $C(p)$, $i(M)$ and $d(M)$ for an arbitrary Riemannian manifold $M$;
but for Riemannian symmetric spaces, the task is much easier.
Richard Crittenden discussed conjugate points and cut points in symmetric spaces in \cite{Cr}; where
he claimed that the conjugate locus is determined by the \textit{diagram} of a single Cartan subalgebra and the isotropy
group, and proved that the cut locus of $p$ coincides with the first conjugate locus of $p$ for every $p\in M$
when $M$ is \textit{simply connected} (Cheeger gave another proof in \cite{Ch}). Based on his work, the author computed $i(M)$ and
$d(M)$ for every type of simply connected, compact and irreducible Riemannian symmetric spaces according to the corresponding
Dynkin diagram and Satake diagram in \cite{Y}.
Now, we study the non-simply connected cases.
The purpose of the paper is to determine cut locus of an arbitrary non-simply connected, compact and irreducible Riemmannian symmetric
space; as an application, we compute $i(M)$ and $d(M)$ for every type of them. The author hopes the results be beneficial
to doing further research for general geometric properties on symmetric spaces of compact type.

In Section 2, we summarize the results about cut locus of an arbitrary simply connected, compact Riemannian symmetric space,
which are due to Richard Crittenden; but our denotation is mainly from \cite{Hel} and \cite{AB}.
Notice the concept of \textit{Cartan polyhedron}, which
plays an important role in the expression of the cut locus and the computation of $i(M)$ and $d(M)$ for both simply connected case
(cf. \cite{Y}) and non-simply connected case. Moreover, we compute the kernel of the exponential mapping explicitly and give two easily-seen corollaries, which are
useful for the next sections.

E.Cartan and M.Takeuchi have studied the fundamental group of compact Riemannian symmetric spaces,
see \cite{Ta}. But for the expression of the cut locus, we adopt a new idea of describing the fundamental group. At the beginning of
Section 3, we explore the relationship between $Z_{\wtd{M}}(\wtd{K})$ and the restricted root system,
where $\wtd{M}=\wtd{U}/\wtd{K}$ is the universal covering space of $M$ and $Z_{\wtd{M}}(\wtd{K})$ denotes
the points in $\wtd{M}$ invariant under the left action of $\wtd{K}$;  then we claim that
there is a one-to-one correspondence between every subgroup of $Z_{\wtd{M}}(\wtd{K})$
and every symmetric space which is locally isometric to $\wtd{M}$, whose fundamental group is isomorphic to
the corresponding subgroup of $Z_{\wtd{M}}(\wtd{K})$.

Then in Section 4, we bring in new denotation (i.e., $P_\G$ and $P'_\G$, where $\G$ is an arbitrary subgroup
of $Z_{\wtd{M}}(\wtd{K})$) and obtain Theorem 4.1 about cut locus, the main Theorem in the paper.

Section 5-8 is the process of computing $i(M)$ and $d(M)$. In Section 5, we compute $(e_i,e_j)$ for every type
of $\Si$ (restricted root system), where $e_1,\cdots, e_l$ denote the vertices of Cartan polyhedron, $(,)$
denotes the Killing form; and give
the group structure of $Z_{\wtd{M}}(\wtd{K})$, which is completely determined by $\Si$. In Section 6, we introduce two
new qualities, i.e., $i(P_\G)$ and $d(P_\G)$ and express them in the form of $(\psi,\psi)$, where $\psi$ is
the highest restricted root ; later in Section 7, we compute $(\psi,\psi)$ for every type of reduced, compact and irreducible
orthogonal involutive Lie algebra (the work is first done by X.S. Liu in \cite{L}, and our method is similar to
\cite{Y}, so we omit the details of computation); then combining the results of Section 6 and Section 7, $i(P_\G)$ and
$d(P_\G)$ are determined explicitly. In Section 8, we give the geometric meaning of a parameter $\ep>0$, which only depend on the metric of
$M$, and then we list $i(M)$ and $d(M)$ for every type
of non-simply connected, compact and irreducible Riemannian symmetric space when $\ep=1$, $Ric=1/2$ in Table 8.1 and Table 8.2 on the basis
of what we have done in Section 5-7. However, when $\wtd{M}=SU(n)/SO(n), SU(2n)/Sp(n)$ or $SU(n)$, $\G=\Z_p$
such that $2<p<n$, the author temporarily have no idea to compute $d(\wtd{M}/\G)$. Our computation is
on the basis of the Dynkin diagram of every reduced root system and the \textit{Satake diagram}
of every type of reduced, compact and irreducible
orthogonal involutive Lie algebra given by Araki in \cite{Ar}.

\section{Some results about the cut locus }
\setcounter{equation}{0}

Let $\sw{u}$ be a compact semisimple Lie algebra and $\th$ an involutive automorphism of $\sw{u}$, then $\th$ extends
uniquely to a complex involutive automorphism of $\sw{g}$, the complexification of $\sw{u}$. We have then the direct decompositions
\begin{eqnarray}
\sw{u}=\sw{k}_0\oplus \sw{p}_*;\qquad \mbox{where }\sw{k}_0=\{X\in \sw{u}:\th(X)=X\},\sw{p}_*=\{X\in \sw{u}:\th(X)=-X\}.
\end{eqnarray}
Let $\lan,\ran$ be an inner product on $\sw{p}_*$ invariant under $Ad\ \sw{k}_0$, then $(\sw{u},\th,\lan,\ran)$
is an orthogonal involutive Lie algebra; without loss of generality we can assume it is \textit{reduced} (cf. \cite{AB} pp.20-21).
Let $M=U/K$ with $U$-invariant metric $g$ is a compact
Riemannian symmetric space which associates with $(\sw{u},\th,\lan,\ran)$, then there is a natural correspondence
between $(T_o M,g)$ and $(\sw{p}_*,\lan,\ran)$, where $o=eK$; in the following text we identify $T_o M$ and $\sw{p}_*$.

It is well known that the geodesic emanating from $o$ with tangent vector $X\in \sw{p}_*$ is given by $\g(t)=\exp(tX)K$,
where $t\ra \exp(tX)$ is a one-parameter subgroup of $U$ (see \cite{Hel} p.208); i.e., if we denote by $\Exp:\sw{p}_*\ra (M,g)$ the exponential mapping, then
$\Exp(X)=\exp(X)K$; and
\begin{eqnarray}
d\ \Exp_X=d\tau(\exp X)_o\circ \sum_{n=0}^\infty \f{(T_X)^n}{(2n+1)!}\qquad X\in \sw{p}_*;
\end{eqnarray}
where $\tau(a)$ denotes the mapping $bK\mapsto abK$ of $U/K$ onto itself for arbitrary $a\in U$, $T_X$ denotes
the restriction of $(ad\ X)^2$ to $\sw{p}_*$ (see \cite{Hel} p. 215). By the properties of compact Lie algebra, $ad\ X$
is anti-symmetric with respect to $\lan,\ran$, thus $T_X$ is symmetric with respect to $\lan,\ran_{\sw{p}_*}$; which yields
that the eigenvalues of $T_X$ are all real; denote by $(\sw{p}_*)_\la(T_X)$ the eigenspace associated to the eigenvalue
$\la$ of $T_X$, then obviously
\begin{eqnarray}
\sum_{n=0}^\infty \f{(T_X)^n}{(2n+1)!}\Big|_{(\sw{p}_*)_\la(T_X)}=\left\{\begin{array}{cc} \mathbf{1} & \la=0;\\ \f{1}{\sqrt{\la}}\sinh(\sqrt{\la})& \la>0;\\ \f{1}{\sqrt{-\la}}\sin(\sqrt{-\la}) & \la<0.\end{array}\right.
\end{eqnarray}
Therefore
\begin{eqnarray}
\ker(d\ \Exp_X)=\bigoplus_{\la<0,\sqrt{-\la}\in \pi\Z}(\sw{p}_*)_\la(T_X).
\end{eqnarray}

Let $\sw{h}_{\sw{p}_*}$ denote an arbitrary maximal abelian subspace of $\sw{p}_*$, $\sw{h}_{\sw{k}_0}$ be an abelian subalgebra
of $\sw{k}_0$ such that $\sw{h}_{\sw{k}_0}\oplus \sw{h}_{\sw{p}_*}$ is a maximal abelian subalgebra of $\sw{u}$,
and $\sw{h}$ denote the subalgebra of $\sw{g}$ generated by $\sw{h}_{\sw{k}_0}\oplus \sw{h}_{\sw{p}_*}$. Denote
$\sw{p}_0=\sqrt{-1}\sw{p}_*$, $\sw{p}=\sw{p}_*\otimes \C$, $\sw{k}=\sw{k}_0\otimes \C$, $\sw{h}_{\sw{p}_0}=\sqrt{-1}\sw{h}_{\sw{p}_*}$,
$\sw{h}_{\sw{p}}=\sw{h}_{\sw{p}_*}\otimes \C$, $\sw{h}_{\sw{k}}=\sw{h}_{\sw{k}_0}\otimes \C$,
then the Killing form $(,)=B(,)$ is positive on $\sqrt{-1}\sw{h}_{\sw{k}_0}\oplus \sw{h}_{\sw{p}_0}$;
let $\De$ be the root system of $\sw{g}$ with respect to $\sw{h}$, then $\sqrt{-1}\sw{h}_{\sw{k}_0}\oplus \sw{h}_{\sw{p}_0}$
is the real linear space generated by $\De$, which is denoted by $\sw{h}_\R$. Denote by $\De^+$ the subset of $\De$ formed by the positive roots with
respect to a lexicographic ordering of $\De$; for every $\a\in \De$, denote by $\a^\th=\th(\a)$, by $\bar{\a}=1/2(\a-\a^\th)$
the orthogonal projection of $\a$ into $\sw{p}_0$. Denote by $\De_0=\{\a\in \De:\bar{\a}=0\}$, $\De_{\sw{p}}=\{\a\in \De:\bar{\a}\neq 0\}$, $P_+=\De^+\cap
\De_{\sw{p}}$; by $\Si=\{\bar{\a}:\a\in \De_{\sw{p}}\}$ the \textit{restricted root system}. $\Si$
has a compatible ordering with $\De$, and $\Si^+=\{\bar{\a}:\a\in P_+\}$. Denote by
\begin{eqnarray}
&&\sw{g}_\g=\{x\in \sw{g}:[H,x]=(H,\g)x,H\in \sw{p}\}\qquad \g\in \Si,\\
&&\sw{k}_\g=(\sw{g}_\g\oplus \sw{g}_{-\g})\cap \sw{k},\ \sw{p}_\g=(\sw{g}_\g\oplus \sw{g}_{-\g})\cap \sw{p}\qquad \g\in \Si^+,
\end{eqnarray}
and by $m_\g=\dim_\C \sw{g}_\g$ the \textit{multiplicity} of $\g$, then
\begin{eqnarray}
\sw{g}=\sw{z}_{\sw{g}}(\sw{h}_\sw{p})\oplus \Big(\bigoplus_{\g\in \Si}\sw{g}_\g\Big),\
\sw{k}=\sw{z}_{\sw{k}}(\sw{h}_\sw{p})\oplus \Big(\bigoplus_{\g\in \Si^+}\sw{k}_\g\Big),\
\sw{p}=\sw{h}_{\sw{p}}\oplus \Big(\bigoplus_{\g\in \Si^+}\sw{p}_\g\Big)
\end{eqnarray}
and
\begin{eqnarray}
m_\g=\big|\{\a\in \De_{\sw{p}}:\bar{\a}=\g\}\big|,\ \dim \sw{k}_\g=\dim \sw{p}_\g=m_\g
\end{eqnarray}
(cf. \cite{Hel} pp.283-293).

For every $X\in \sw{p}_*$, there exists $k\in K$ and $H\in \sw{h}_{\sw{p}_*}$, such that $X=Ad(k)H$ (cf. \cite{AB} p. 31).
For arbitrary $u\in \sw{p}_\g$, $(ad\ H)^2 u=-\big(ad\ (-\sqrt{-1}H)\big)^2 u=-(-\sqrt{-1}H,\g)^2 u$;
$-\sqrt{-1}H,\g\in \sw{h}_{\sw{p}_0}$ yields $(-\sqrt{-1}H,\g)^2\geq 0$; i.e.,
\begin{eqnarray}
\sw{p}_\g\cap \sw{p}_*\subset (\sw{p}_*)_{-(-\sqrt{-1}H,\g)^2}(T_H).
\end{eqnarray}
Since $X=Ad(k)H$, $ad\ X=Ad(k)\circ ad\ H\circ Ad(k)^{-1}$ and moreover
$T_X=(ad\ X)^2=Ad(k)\circ (ad\ H)^2\circ Ad(k)^{-1}=Ad(k)\circ T_H\circ Ad(k)^{-1}$; which yields the eigenvalues of
$T_X$ coincide with the eigenvalues of $T_H$ and for every eigenvalue $\la$, $(\sw{p}_*)_\la(T_X)=Ad(k)\big((\sw{p}_*)_\la(T_H)\big)$.
By (2.4), (2.7), (2.9), we have

\begin{Thm}
Let $M=U/K$ be a compact Riemannian symmetric space such that $U$ is a semi-simple and compact Lie group,
and the denotation of $\sw{p}_*,\sw{k}_0,\sw{h}_{p_*},\Si, \Exp$ is similar to above, then
for every $X=Ad(k)H\in \sw{p}_*$, where $k\in K,H\in \sw{h}_{p_*}$, $X$ is a conjugate point in $T_o M$
if and only if there exists at least one $\g\in \Si$, such that
\begin{eqnarray}
(H,\g)\in \pi\sqrt{-1}(\Z-0)
\end{eqnarray}
and $\ker(\Exp)_X$ is the direct sum of $Ad(k)(\sw{p}_\g\cap \sw{p}_*)$ such that $\g\in \Si^+$ and $(H,\g)\in \pi\sqrt{-1}(\Z-0)$.
\end{Thm}

Denote by $C(p)$ the cut locus of $p\in M$ in $T_p M$, by
\begin{eqnarray}
\mathfrak{S}(p)=\{X\in T_p M: d\big(p,\exp_p(X)\big)=|X|\};
\end{eqnarray}
 then $X\in \mathfrak{S}(p)$ if and only if there exists $X_0\in C(p)$ and $t\in [0,1]$ such that
$X=tX_0$, and moreover $C(p)=\p \mathfrak{S}(p)$ (cf. \cite{CE} pp. 94-95). In 1962, Richard Crittenden proved the following proposition in \cite{Cr}:

\begin{Lem}
Let $M$ be a simply connected complete symmetric space, for every $p\in M$, the cut locus of
$p$ coincides with the first conjugate locus of $p$.
\end{Lem}

Then by Lemma 2.1 and Theorem 2.1, $X=Ad(k)H\in \mathfrak{S}(o)$ if and only if $(tH,\g)\notin \pi\sqrt{-1}(\Z-0)$
for every $t\in [0,1)$ and $\g\in \Si$, where $k\in K$ and $H\in \sw{h}_{\sw{p}_*}$; which implies
\begin{eqnarray}
-\pi\sqrt{-1}\leq (H,\g)\leq \pi\sqrt{-1}\qquad\mbox{ for every }\g\in \Si.
\end{eqnarray}

Now we denote by $C$ the Weyl chamber with respect the ordering of $\Si$, i.e., $C=\{x\in \sw{h}_{\sw{p}_0}:
(x,\g)>0\mbox{ for every }\g\in \Si^+\}$, by $\Pi$ the set of simple roots.
Recall that the planes
$(x,\g)\in \Z$($\g\in \Si$) in $\sw{h}_{\sw{p}_0}$ constitute the \textit{diagram} $D(\Si)$ of $\Si$, and the closure of a connected
component of $\sw{h}_{\sw{p}_0}-D(\Si)$ will be called a \textit{Cartan polyhedron}. Especially,
let $A$ be the set of maximal roots, then the inequalities $(x,\g)\geq 0(\g\in \Pi), (x,\be)\leq 1(\be\in A)$ define
a Cartan polyhedron, which is denoted by $\triangle$ (See \cite{AB} p. 10). Obviously $\triangle\subset \ol{C}$, where $\ol{C}$
denotes the closure of $C$ in $\sw{h}_{\sw{p}_0}$. Since
Weyl group $W$ permutes Weyl chamber in a simply transitive manner  and every element of Weyl group
can be extended to $Ad_{\sw{u}}(\sw{k}_0)$ (See \cite{Hel} pp. 288-290), for every $X\in \sw{p}_*$, there exist
$k\in K$ and $H\in \sqrt{-1}\ \ol{C}$ such that $X=Ad(k)H$. By (2.12), $X\in \mathfrak{S}(o)$ if and only if
$H\in \pi\sqrt{-1}\triangle$. Then we have

\begin{Thm}
Let $M=U/K$ be a simply connected and compact Riemannian symmetric space such that $U$ is a semi-simple and compact Lie group,
and the denotation of $\sw{p}_*,\sw{k}_0,\sw{h}_{p_*},\Si,\triangle$ is similar to above, then $\mathfrak{S}(o)=Ad(K)
(\pi\sqrt{-1}\triangle)$.
\end{Thm}

From Theorem 2.2, by the completeness of $M$, we easily obtain the following Corollaries:

\begin{Cor}
The assumption and denotation are similar to Theorem 2.2, then for every $p\in M$, there exists $k\in K$ and
$x\in \triangle$, such that $p=\Exp\big(Ad(k)(\pi\sqrt{-1}x)\big)$, and $d(o,p)=\pi|\sqrt{-1}x|$; where
$|X|=\lan X,X\ran^{1/2}$ for arbitrary $X\in T_o M=\sw{p}_*$.
\end{Cor}

\begin{Cor}
The assumption and denotation are similar to Theorem 2.2 and Corollary 2.1; denote $p=\Exp(\pi\sqrt{-1}x)$,
$q=\Exp(\pi\sqrt{-1}y)$, where $x,y\in \triangle$, then $d(p,q)=\pi|\sqrt{-1}(y-x)|$.
\end{Cor}

Proof. Since the metric $g$ on $M$ is $U$-invariant,
\begin{eqnarray}
d(p,q)=d\big(\tau(\exp(-\pi\sqrt{-1}x))p,\tau(\exp(-\pi\sqrt{-1}x))q\big)=d\big(o,\Exp(\pi\sqrt{-1}(y-x))\big);
\end{eqnarray}
since $x,y\in \triangle$, for every $\g\in \Si^+$, $(x,\g),(y,\g)\in [0,1]$, thus $(y-x,\g)\in [-1,1]$;
then (2.12) yields $\pi\sqrt{-1}(y-x)\in \mathfrak{S}(o)$; by the definition of $\mathfrak{S}(o)$ and (2.13),
$d(p,q)=\pi|\sqrt{-1}(y-x)|$. $\Box$

\section{Some properties of $Z_{\wtd{M}}(\wtd{K})$}
\setcounter{equation}{0}

In this section, we assume $(\sw{u},\th,\lan,\ran)$ be a reduced, compact and irreducible orthogonal involutive Lie
algebra,  $\wtd{U}$ be the simply connected Lie group associated with $\sw{u}$, and
$\wtd{M}=\wtd{U}/\wtd{K}$ with $\wtd{U}$-invariant metric $\wtd{g}$ be a simply connected Riemannian symmetric space
associated with $(\sw{u},\th,\lan,\ran)$. Denote by $\wtd{\si}$ the involutive automorphism of $\wtd{U}$ induced by
$\th$, then the fixed point set of $\wtd{\si}$ (denoted by $\wtd{U}_{\wtd{\si}}$) is connected (see \cite{RH},\cite{Ch},\cite{CE} pp.102-103); which yields that
$\wtd{K}=\wtd{U}_{\wtd{\si}}$ is the connected Lie subgroup of $\wtd{U}$ generated by $\sw{k}_0$. Denote by $\wtd{\exp}$ the exponential
mapping of $\sw{u}$ onto $\wtd{U}$, by $\wtd{\Exp}:\sw{p}_*\ra \wtd{M}$ $X\mapsto \wtd{\exp}(X)\wtd{K}$,
by $\wtd{o}=\wtd{e}\wtd{K}$, where $\wtd{e}$ the identity element of $\wtd{U}$. The denotation
of $\sw{h}_{\sw{p}_*},\sw{h}_{\sw{p}_0},\sw{h}_{\sw{p}},\sw{h}_{\sw{k}_0},\sw{h}_{\sw{k}},\Si,\Pi,\triangle$ is similar to Section 2;
since $(\sw{u},\th,\lan,\ran)$ is irreducible, $\Si$ is also irreducible and $\triangle$ is a simplex; let $\psi$
be \textit{the highest restricted root}, $\Pi=\{\g_1,\cdots,\g_l\}$ ($l=\rank(\Si)=\dim(\sw{h}_{\sw{p}_*})$), and $d_1,\cdots,d_l
\in \Z^+$ such that $\psi=\sum_{i=1}^l d_i\g_i$, then the vertices of $\triangle$ include
\begin{eqnarray}
e_1,\cdots,e_l;\qquad (e_j,\g_i)=\f{1}{d_j}\de_{ij}.
\end{eqnarray}

Denote $Z_{\wtd{M}}(\wtd{K})=\{p\in \wtd{M}:\tau(k)p=p \mbox{  for every } k\in \wtd{K}\}$, then on $Z_{\wtd{M}}(\wtd{K})$
we have the following proposition:

\begin{Pro}
There exists a natural group structure on $Z_{\wtd{M}}(\wtd{K})$ if we define $a\wtd{K}\cdot b\wtd{K}=ab\wtd{K}$
and $(a\wtd{K})^{-1}=a^{-1}\wtd{K}$ for every $a\wtd{K},b\wtd{K}\in Z_{\wtd{M}}(\wtd{K})$. Then $Z_{\wtd{M}}(\wtd{K})$
is a finite abelian group, and for every $p=a\wtd{K}\in \wtd{M}-\{\wtd{o}\}$, the following conditions
are equivalence:

(a) $p\in Z_{\wtd{M}}(\wtd{K})$;

(b) $a\in N_{\wtd{U}}(\wtd{K})$, where $N_{\wtd{U}}(\wtd{K})$ denotes the normalizer of $\wtd{K}$ in $\wtd{U}$;

(c) $aa^*\in Z(\wtd{U})\cap \wtd{\exp}(\sw{h}_{\sw{p}_*})$, where $Z(\wtd{U})$ denotes the center of $\wtd{U}$ and $a^*=\wtd{\si}(a)^{-1}$;

(d) $p=\wtd{\Exp}(\pi\sqrt{-1}e_{j})$ such that $d_{j}=1$.
\end{Pro}

Proof.
$(a)\Longleftrightarrow (b)$: If $p=a\wtd{K}\in Z_{\wtd{M}}(\wtd{K})$, then for every $k\in \wtd{K}$, $a\wtd{K}
=ka\wtd{K}$, which yields $a^{-1}ka\in \wtd{K}$, i.e., $a\in N_{\wtd{U}}(\wtd{K})$; and vice versa.

$(b)\Longrightarrow (c)$: It is well known that $\wtd{U}=\wtd{K}\wtd{\exp}(\sw{h}_{\sw{p}_*})\wtd{K}$ (cf. \cite{AB} pp. 74-76), so
$a=b_1\wtd{\exp}X b_2$ for some $X\in \sw{h}_{\sw{p}_*}$ and $b_1,b_2\in \wtd{K}$.
$a\in N_{\wtd{U}}(\wtd{K})$ yields $\wtd{\exp}X\in N_{\wtd{U}}(\wtd{K})$; by the easily-seen facts that $(bc)^*=c^*b^*$
and $k^*=k^{-1}$ for arbitrary $b,c\in \wtd{U}$ and $k\in \wtd{K}$, we have
$$aa^*=b_1\wtd{\exp}Xb_2{b_2}^*(\wtd{\exp}X)^*b_1^*=F_{b_1}\big(\wtd{\exp}X(\wtd{\exp} X)^*\big);\qquad (\mbox{where }F_b(c)=bcb^{-1})$$
so $aa^*\in Z(\wtd{U})\cap \wtd{\exp}(\sw{h}_{\sw{p}_*})$ if and only if $\wtd{\exp}X(\wtd{\exp}X)^*$ does;
without loss of generality we can assume $a=\wtd{\exp}X$. For every $k\in \wtd{K}$, there exists $k'\in \wtd{K}$,
such that $k a=ak'$, thus
$$F_k(aa^*)=kaa^*k^{-1}=ka(ka)^*=(ak')(ak')^*=aa^*.$$
Since $\sw{h}_{\sw{p}_*}$ is abelian, $aa^*=\wtd{\exp}(2X)$ is invariant under $F_{\wtd{\exp}Y}$ for arbitrary $Y\in \sw{h}_{\sw{p}_*}$; furthermore,
it is invariant under $F_b=F_{k_1}\circ F_{\wtd{\exp}Y}\circ F_{k_2}$ for arbitrary $b=k_1\wtd{\exp}Yk_2\in \wtd{U}$. Hence (c) holds.

$(c)\Longrightarrow (b)$: Denote $a=b_1\wtd{\exp}X b_2$, where $b_1,b_2\in \wtd{K}$ and $X\in \sw{h}_{\sw{p}_*}$. $Z(\wtd{U}) \ni aa^*
=b_1\wtd{\exp}(2X)b_1^{-1}$ implies $\wtd{\exp}(2X),\wtd{\exp}(-2X)\in Z(\wtd{U})$; then for every $k\in \wtd{K}$
$$\wtd{\si}(F_{\wtd{\exp}X}k)=F_{\wtd{\si}(\wtd{\exp}X)}\wtd{\si}(k)=F_{\wtd{\exp}(-X)}k=F_{\wtd{\exp}(-2X)}F_{\wtd{\exp}X}k=F_{\wtd{\exp}X}k;$$
i.e., $F_{\wtd{\exp}X}k\in \wtd{K}$, $\wtd{\exp}X\in N_{\wtd{U}}(\wtd{K})$. Hence $a=b_1\wtd{\exp}(X)b_2\in N_{\wtd{U}}(\wtd{K})$.

$(a)\Longrightarrow (d)$: By Corollary 2.1, there exists $k\in \wtd{K}$ and $x\in \triangle$, such that
$p=\wtd{\Exp}\big(Ad(k)(\pi\sqrt{-1}x)\big)\\=\tau(k)\wtd{\Exp}(\pi\sqrt{-1}x)$, then (a) implies $p=\tau(k^{-1})p
=\wtd{\Exp}(\pi\sqrt{-1}x)$.

Denote $X=\pi\sqrt{-1}x$, then for arbitrary $Y\in \sw{k}_0$ and $t\in \R$,
\begin{eqnarray}
\wtd{\Exp}(X)=p=\tau\big(\wtd{\exp}(t Y)\big)p=\tau\big(\wtd{\exp}(t Y)\big)\wtd{\Exp}(X)
=\wtd{\Exp}\big(\exp(t\ ad Y)X\big).
\end{eqnarray}
Differentiate both sides of (3.2) and then let $t=0$, we have
\begin{eqnarray}
(d\ \wtd{\Exp})_X[Y,X]=0.
\end{eqnarray}
Notice that $x\in \triangle$; applying Theorem 2.1, we obtain
\begin{eqnarray}
[Y,X]\in \bigoplus_{\g\in \Si^+,(x,\g)=1}(\sw{p}_\g\cap \sw{p}_*).
\end{eqnarray}
$p\neq \wtd{o}$ yields $X\neq 0$, then there exists $\g_j\in \Pi$, $(X,\g_j)\neq 0$; take nonzero $Y\in \sw{k}_{\g_j}\cap \sw{k}_0$,
then $[Y,X]\in \sw{p}_{\g_j}\cap \sw{p}_*$ and (3.4) yields $(x,\g_j)=1$; since $x\in \triangle$, we have $(x,\psi)=1$
and moreover
$$1=(x,\psi)=\sum_{i=1}^l d_i(x,\g_i);$$
which yields $d_j=1$ and $(x,\g_i)=\de_{ij}$; i.e., $x=e_j$.

$(d)\Longrightarrow (c)$: $p=\wtd{\Exp}(\pi\sqrt{-1}e_j)$ yields $aa^*=\wtd{\exp}(2\pi\sqrt{-1}e_j)\in \wtd{\exp}(\sw{h}_{\sw{p}_*})$.
Denote by $Ad:\wtd{U}\ra GL(\sw{u})$ the adjoint homomorphism, then
\begin{eqnarray}
&&Ad\big(aa^{*}\big)|_{\sw{g}_\g}=\exp(2\pi\sqrt{-1}ad\ e_j)|_{\sw{g}_\g}=e^{2\pi\sqrt{-1}(e_j,\g)}=\mathbf{1}\qquad \g\in \Si;\nonumber\\
&&Ad\big(aa^{*}\big)|_{\sw{z}_{\sw{g}}(\sw{h}_\sw{p})}=\exp(2\pi\sqrt{-1}ad\ e_j)|_{\sw{z}_{\sw{g}}(\sw{h}_\sw{p})}=\mathbf{1}.
\end{eqnarray}
By (2.7), $aa^{*}\in \ker(Ad)\subset Z(\wtd{U})$.
\newline

$N_{\wtd{U}}(\wtd{K})$ is a Lie subgroup of $\wtd{U}$, and the Lie algebra associated to $N_{\wtd{U}}(\wtd{K})$ is
$\sw{n}_{\sw{u}}(\sw{k}_0)=\sw{k}_0$ (since $(\sw{u},\th,\lan,\ran)$ is semi-simple, cf. \cite{AB} p.25); so $\wtd{K}$ is the identity
component of $N_{\wtd{U}}(\wtd{K})$ and then $Z_{\wtd{M}}(\wtd{K})=N_{\wtd{U}}(\wtd{K})/\wtd{K}$ is a finite group.

Define $\Psi:Z_{\wtd{M}}(\wtd{K})\ra Z(\wtd{U})\cap \wtd{\exp}(\sw{h}_{\sw{p}_*})$
\begin{eqnarray}
a\wtd{K}\mapsto aa^*;
\end{eqnarray}
obviously
$\Psi(\wtd{o})=\wtd{e}$ and
$$\Psi(a\wtd{K}\cdot b\wtd{K})=\Psi(ab\wtd{K})=ab(ab)^*=a(bb^*)a^*=aa^*(bb^*)=\Psi(a\wtd{K})\Psi(b\wtd{K})$$
for every $a\wtd{K},b\wtd{K}\in Z_{\wtd{M}}(\wtd{K})$; if $\Psi(a\wtd{K})=\wtd{e}$, then $\wtd{\si}(a)=a$
and therefore $a\wtd{K}=\wtd{o}$; hence $\Psi$ is a monomorphism. $Z_{\wtd{M}}(\wtd{K})$ could be considered a subgroup of $Z(\wtd{U})$, which is an abelian group.
$\Box$
\newline

By Corollary 2.2, $d\big(\wtd{\Exp}(\pi\sqrt{-1}e_j),\wtd{\Exp}(\pi\sqrt{-1}e_k)\big)=\pi|\sqrt{-1}(e_k-e_j)|\neq 0$
when $j\neq k$, then Proposition 3.1 tells us
\begin{eqnarray}
Z_{\wtd{M}}(\wtd{K})=\big\{\wtd{\Exp}(\pi\sqrt{-1}e_j):d_j=1\big\}\cup \{\wtd{o}\}
\end{eqnarray}
For every $j\neq k$ such that $d_j=d_k=1$, $(e_j-e_k,\g)\in [-1,1]$ for arbitrary $\g\in \Si$, so there exists $\om\in W$
and $x\in \triangle$, such that $e_j-e_k=\om(x)$; let $k\in \wtd{K}$ such that $\om=Ad(k)|_{\sw{h}_{\sw{p}_0}}$,
then
\begin{eqnarray*}
Z_{\wtd{M}}(\wtd{K})&\ni& \wtd{\Exp}(\pi\sqrt{-1}e_j)\wtd{\Exp}(\pi\sqrt{-1}e_k)^{-1}=\wtd{\Exp}\big(\pi\sqrt{-1}(e_j-e_k)\big)\\
&=&\tau(k)\wtd{\Exp}(\pi\sqrt{-1}x)=\wtd{\Exp}(\pi\sqrt{-1}x).
\end{eqnarray*}
From the proceed of proving Proposition 3.1, we have $x=e_r$ for some $1\leq r\leq l$ such that $d_r=1$. Similarly,
there exists $1\leq s\leq l$ such that $d_s=1$ and $\wtd{\Exp}(\pi\sqrt{-1}e_j)^{-1}=\wtd{\Exp}(\pi\sqrt{-1}e_s)$.
It tells us that the group structure of $Z_{\wtd{M}}(\wtd{K})$ can be uniquely determined by the type of $\Si$,
and using the technology stated above, we can write $\wtd{\Exp}(\pi\sqrt{-1}e_j)\wtd{\Exp}(\pi\sqrt{-1}e_k)^{-1}$
and $\wtd{\Exp}(\pi\sqrt{-1}e_j)^{-1}$ precisely. We shall give the results for every type of
$\Si$ in Section 5 after concrete computation.

$Z_{\wtd{M}}(\wtd{K})$ has a close relationship with the fundamental group of Riemannian symmetric spaces.
Let $M=U/K$ with $U$-invariant metric be a non-simply connected Riemannian symmetric space associated with
$(\sw{u},\th,\lan,\ran)$, then the universal covering group of $U$ is $\wtd{U}$; denote by $\chi:\wtd{U}\ra U$
the covering homomorphism and by $\pi:\wtd{M}=\wtd{U}/\wtd{K}\ra M=U/K$
\begin{eqnarray}
a\wtd{K}\ra \chi(a)K;
\end{eqnarray}
then $\pi$ is obviously a covering map, and the pullback metric $\pi^{-1}g$ coincides with $\wtd{g}$. In this
case, $M$ is called a \textit{Clifford-Klein form} of $\wtd{M}$ ; $M$ is isomorphic to the quotient
of $\wtd{M}$ by a properly discontinuous group of isometries $L$, which is isomorphic to $\pi^{-1}(o)=\chi^{-1}(K)/\wtd{K}$
(cf. \cite{AB} pp.101-105).
By $\chi(\wtd{K})\subset K$, $\chi^{-1}(K)/\wtd{K}$ is invariant under $\tau(k)$ for arbitrary $k\in \wtd{K}$;
furthermore, since $\chi^{-1}(K)/\wtd{K}$ is discrete, any point of which is invariant under $\tau(k)$; therefore
$\pi^{-1}(o)\subset Z_{\wtd{M}}(\wtd{K})$, i.e., the fundamental group of $M$ is a subgroup of $Z_{\wtd{M}}(\wtd{K})$.

Conversely, let $\G$ be an arbitrary subgroup of $Z_{\wtd{M}}(\wtd{K})$; by Proposition 3.1, $\Psi(\G)$ is a subgroup of
$Z(\wtd{U})$, where $\Psi$ is defined in (3.6); let $U=\wtd{U}/\Psi(\G)$, denote by $\chi:\wtd{U}\ra U$ the covering
homomorphism. For any $aa^*\in \Psi(\G)$ such that $a\wtd{K}\in \G$, since $\wtd{\exp}(b)=(b^*)^{-1}$ and $b^{**}=b$,
$$\wtd{\si}(aa^*)=\big((aa^*)^*\big)^{-1}=(a^{**}a^*)^{-1}=(aa^*)^{-1}\in \Psi(\G),$$
i.e., $\wtd{\si}$ keeps $\Psi(\G)$ invariant; so $\wtd{\si}$ induces a involutive automorphism of $U$,
which is denoted by $\si$. Let $K=U_\si$, and the definition of $\pi$ is similar to (3.13), then for every $a\in \wtd{U}$,
$a\wtd{K}\in \pi^{-1}(o)=\chi^{-1}(K)/\wtd{K}$ if and only if
$$\si\big(\chi(a)\big)=\chi(a)\qquad \mbox{i.e., } \Psi(a\wtd{K})=aa^*=a\wtd{\si}(a)^{-1}\in \ker(\chi)=\Psi(\G);$$
so the fundamental group of $M$ is isomorphic to $\G$. We can express $M$ as $\wtd{M}/\G$.

Therefore,
all of the subgroups of $Z_{\wtd{M}}(\wtd{K})$, which is uniquely determined by the type of $\Si$, could
completely determine every  compact and irreducible Riemannian symmetric space which is locally isometric to
$\wtd{M}$, i.e., every Clifford-Klein form of $\wtd{M}$.

\section{The cut locus of non-simply connected, compact and irreducible Riemannian symmetric spaces}
\setcounter{equation}{0}

Our assumption and the denotation of $(\sw{u},\th,\lan,\ran)$, $\wtd{M}$, $Z_{\wtd{M}}(\wtd{K})$ is similar to Section 3.
Let $M=\wtd{M}/\G$ be a Clifford-Klein form of $\wtd{M}$,
where $\G$ is a subgroup of $Z_{\wtd{M}}(\wtd{K})$ satisfying $\G\neq \{\wtd{o}\}$; and denote by $\pi:\wtd{M}\ra M$
the covering mapping. We will discuss the cut locus of $M$; and our denotation of $C(p)$ and $\mathfrak{S}(p)$ is similar to Section 2.

Obviously $\pi\big(\wtd{\Exp}(X)\big)=\Exp(X)$ for arbitrary $X\in \sw{p}_*$. By the properties
of covering maps, we have
\begin{eqnarray}
d_M\big(o,\Exp(X)\big)= \min_{p\in \G}d_{\wtd{M}}\big(p,\wtd{\Exp}(X)\big).
\end{eqnarray}
Let $X\in \mathfrak{S}(o)$, i.e., $d_M\big(o,\Exp(X)\big)=|X|$, then (4.1) implies $d_{\wtd{M}}\big(\wtd{o},\wtd{\Exp}(X)\big)\geq |X|$;
on the other hand, obviously $d_{\wtd{M}}\big(\wtd{o},\wtd{\Exp}(X)\big)\leq |X|$; then $d_{\wtd{M}}\big(\wtd{o},\wtd{\Exp}(X)\big)=|X|$,
i.e., $X\in \mathfrak{S}(\wtd{o})$. By Theorem 2.2, there exists $k\in \wtd{K}$ and $x\in \triangle$, such that
$X=Ad(k)(\pi\sqrt{-1}x)$. For any $p=\wtd{\Exp}(\pi\sqrt{-1}e_i)\in \G$, by Corollary 2.2,
\begin{eqnarray}
d_{\wtd{M}}\big(p,\wtd{\Exp}(X)\big)&=&d_{\wtd{M}}\big(p,\tau(k)\wtd{\Exp}(\pi\sqrt{-1}x)\big)=d_{\wtd{M}}\big(\tau(k^{-1})p,\wtd{\Exp}(\pi\sqrt{-1}x)\big)\nonumber\\
&=&d_{\wtd{M}}\big(\wtd{\Exp}(\pi\sqrt{-1}e_i),\wtd{\Exp}(\pi\sqrt{-1}x)\big)=\pi\big|\sqrt{-1}(x-e_i)\big|.
\end{eqnarray}
Then by (4.1),
\begin{eqnarray}
&&|X|= \min_{p\in \G}d_{\wtd{M}}\big(p,\wtd{\Exp}(X)\big)
=\min\big\{|X|,\pi\big|\sqrt{-1}(x-e_i)\big|:\wtd{\Exp}(\pi\sqrt{-1}e_i)\in \G\big\}\nonumber\\
&&\ \ \ \ \ =\min\big\{\pi\big|\sqrt{-1}x\big|,\pi\big|\sqrt{-1}(x-e_i)\big|:\wtd{\Exp}(\pi\sqrt{-1}e_i)\in \G\big\};\nonumber\\
&&\mbox{i.e., }x\in \triangle\mbox{ and }\big|\sqrt{-1}x\big|\leq \big|\sqrt{-1}(x-e_i)\big|\mbox{ for every }\wtd{\Exp}(\pi\sqrt{-1}e_i)\in \G.
\end{eqnarray}
Conversely, if (4.3) is satisfied, it is easy to check that $d_M\big(o,\Exp(X)\big)=|X|$. Therefore, (4.3) is a
necessary and sufficient condition for $X\in \mathfrak{S}(o)$.

The condition (4.3) can be simplified further. Since $(\sw{u},\th,\lan,\ran)$ is irreducible, there exists a positive constant
$\ep$, such that
\begin{eqnarray}
\lan,\ran=-\ep(,) \qquad (,)\mbox{ be the Killing form on }\sw{g}
\end{eqnarray}
(cf. \cite{AB} pp.23-26).
Then for every $y\in \sw{p}_0$, $\big|\sqrt{-1}y\big|^2=\ep(y,y)$; hence $\big|\sqrt{-1}x\big|\leq \big|\sqrt{-1}(x-e_i)\big|$
if and only if $\ep(x,x)\leq \ep(x-e_i,x-e_i)$, i.e., $(x,e_i)\leq 1/2(e_i,e_i)$.

As a matter of convenience, we bring in new denotation:

DENOTATION 4.1.
Given an arbitrary subgroup $\G\subset Z_{\wtd{M}}(\wtd{K})$, we denote
\begin{eqnarray}
&&P_\G=\big\{x\in \triangle: (x,e_i)\leq \f{1}{2}(e_i,e_i)\mbox{ for every }\wtd{\Exp}(\pi\sqrt{-1}e_i)\in \G\big\};\\
&&P'_\G=\big\{x\in P_\G: (x,\psi)=1\mbox{ or } (x,e_j)=\f{1}{2}(e_j,e_j)\mbox{ for some }j\mbox{ such that }\wtd{\Exp}(\pi\sqrt{-1}e_j)\in \G\big\}.\nonumber\\
\end{eqnarray}

For any $x\in P_\G$, if $tx\notin P_\G$ for every $t>1$, then $x\in P'_\G$; and vice versa. So we have the following
Theorem:

\begin{Thm}

Let $(\sw{u},\th,\lan,\ran)$ be a reduced, compact and irreducible orthogonal involutive Lie algebra, $\wtd{M}=\wtd{U}
/\wtd{K}$ be the simply connected Riemannian symmetric space associated with $(\sw{u},\th,\lan,\ran)$,
$M=\wtd{M}/\G$ be a Clifford-Klein form of $\wtd{M}$, where $\G$ is a subgroup of $Z_{\wtd{M}}(\wtd{K})$
satisfying $\G\neq \{\wtd{o}\}$, then
$\mathfrak{S}(o)=Ad(K)(\pi\sqrt{-1}P_\G)$ and $C(o)=\p\mathfrak{S}(o)=Ad(K)(\pi\sqrt{-1}P'_\G)$.

\end{Thm}

\section{Some computation on $e_i$ and several corollaries}
\setcounter{equation}{0}

The section is preparation for the next section. One of our purposes is to compute $(e_i,e_j)$,
after which, we will give the group structure of $Z_{\wtd{M}}(\wtd{K})$.

When computing $(e_i,e_j)$, we assume $Z_{\wtd{M}}(\wtd{K})\neq \{\wtd{o}\}$,  which implies $\Si=\sw{a}_l,\sw{b}_l ,
\sw{c}_l,\sw{d}_l,\sw{e}_6$ or $\sw{e}_7$ (by (3.7)). If $\Si$ is a classical root system ($\sw{a}_l$, $\sw{b}_l$, $\sw{c}_l$ or $\sw{d}_l$), then $\Si$
can be imbedded into Euclidean space in a natural manner (see \cite{Hel} pp.461-465); so we can express $e_i$ explicitly according to (3.1)
and then compute $(e_i,e_j)$. Otherwise, when $\Si$ is an exceptional one, the following Lemma takes effect.

\begin{Lem}
The denotation of $\Si,\Pi,\g_i,e_i$ is similar to Section 3, if we denote $\Om_{ij}=(\g_i,\g_j)$, then
$(e_i,e_j)=1/(d_jd_i)(\Om^{-1})_{ji}$.
\end{Lem}

Proof. Since $\{\g_1,\cdots,\g_l\}$ is a basis of $\sw{h}_{\sw{p}_0}$, we can write $e_j=\g_k A_j^k$; then (3.1) yields
\begin{eqnarray}
\f{1}{d_j}\de_{ij}=(e_j,\g_i)=(\g_k A_j^k,\g_i)=A_j^k\Om_{ki};
\end{eqnarray}
so $A_j^k=1/d_j(\Om^{-1})_{ki}\de_{ij}=1/d_j(\Om^{-1})_{kj}$ and
\begin{eqnarray*}
(e_i,e_j)=(\g_k A_i^k,e_j)=\f{1}{d_j}\de_{jk}A_i^k=\f{1}{d_j}\de_{jk}\f{1}{d_i}(\Om^{-1})_{ki}=\f{1}{d_jd_i}(\Om^{-1})_{ji}.\qquad \Box
\end{eqnarray*}

Now we give the detail of computation of $(e_i,e_j)$ for every type of $\Si$.
\newline

$\Si=\sw{a}_l$($l\geq 1$): The corresponding Dynkin diagram is
\setlength{\unitlength}{1mm}
\begin{center}
\begin{picture}(33,5)
\put(1.5,3.5){\circle{2}}
\put(2.5,3.5){\line(1,0){5}}
\put(8.5,3.5){\circle{2}}
\put(9.5,3.5){\line(1,0){5}}
\put(14.8,2.6){$\cdots$}
\put(19,3.5){\line(1,0){5}}
\put(25,3.5){\circle{2}}
\put(26,3.5){\line(1,0){5}}
\put(32,3.5){\circle{2}}
\put(0,0){$\g_1$}
\put(7,0){$\g_2$}
\put(21,0){$\g_{l-1}$}
\put(30.5,0){$\g_l$}
\end{picture}
\end{center}
Denote $\g_i=x_i-x_{i+1}$ ($1\leq i\leq l$), then $\Si=\{\pm(x_i-x_j):1\leq i<j\leq l+1\}$, $\psi=x_1-x_{l+1}=\sum_{i=1}
^l \g_i$ and therefore
$(x_i,x_j)=1/2(\psi,\psi)\de_{ij}$;
by (3.1), we obtain
\begin{eqnarray}
&&e_j=\f{2}{(\psi,\psi)(l+1)}\big((l+1-j)\sum_{k=1}^j x_k-j\sum_{k=j+1}^{l+1}x_k\big)\qquad 1\leq j\leq l;\\
&&(e_i,e_j)=\f{2i(l+1-j)}{(\psi,\psi)(l+1)}\qquad 1\leq i\leq j\leq l.
\end{eqnarray}

$\Si=\sw{b}_l$($l\geq 2$): The corresponding Dynkin diagram is
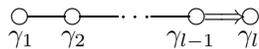
\setlength{\unitlength}{1mm}
\begin{center}
\begin{picture}(34,5)
\put(1.5,3.5){\circle{2}}
\put(2.5,3.5){\line(1,0){5}}
\put(8.5,3.5){\circle{2}}
\put(9.5,3.5){\line(1,0){5}}
\put(14.8,2.6){$\cdots$}
\put(19,3.5){\line(1,0){5}}
\put(25,3.5){\circle{2}}
\put(25.9,2.6){$\Longrightarrow$}
\put(32.2,3.5){\circle{2}}
\put(0,0){$\g_1$}
\put(7,0){$\g_2$}
\put(21,0){$\g_{l-1}$}
\put(30.5,0){$\g_l$}
\end{picture}
\end{center}
Denote $\g_i=x_i-x_{i+1}$($1\leq i\leq l-1$), $\g_l=x_l$, then $\Si=\{\pm(x_i\pm x_j):1\leq i<j\leq l\}\cup
\{\pm x_i:1\leq i\leq l\}$, $\psi=x_1+x_2=\g_1+2\sum_{i=2}^l \g_i$ and therefore
$(x_i,x_j)=1/2(\psi,\psi)\de_{ij}$;
by (3.1), we obtain
\begin{eqnarray}
&&e_1=\f{2}{(\psi,\psi)}x_1,\ e_j=\f{1}{(\psi,\psi)}\sum_{k=1}^j x_k\ (2\leq j\leq l);\\
&&(e_1,e_1)=\f{2}{(\psi,\psi)},\ (e_1,e_j)=\f{1}{(\psi,\psi)}\ (2\leq j\leq l),\ (e_i,e_j)=\f{i}{2(\psi,\psi)}\ (2\leq i\leq j\leq l).\nonumber\\
\end{eqnarray}

$\Si=\sw{c}_l$($l\geq 3$): The corresponding Dynkin diagram is
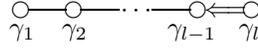
\setlength{\unitlength}{1mm}
\begin{center}
\begin{picture}(34,5)
\put(1.5,3.5){\circle{2}}
\put(2.5,3.5){\line(1,0){5}}
\put(8.5,3.5){\circle{2}}
\put(9.5,3.5){\line(1,0){5}}
\put(14.8,2.6){$\cdots$}
\put(19,3.5){\line(1,0){5}}
\put(25,3.5){\circle{2}}
\put(25.9,2.6){$\Longleftarrow$}
\put(32.2,3.5){\circle{2}}
\put(0,0){$\g_1$}
\put(7,0){$\g_2$}
\put(21,0){$\g_{l-1}$}
\put(30.5,0){$\g_l$}
\end{picture}
\end{center}
Denote $\g_i=x_i-x_{i+1}$($1\leq i\leq l-1$), $\g_l=2x_l$, then $\Si=\{\pm(x_i\pm x_j):1\leq i<j\leq l\}\cup
\{\pm 2x_i:1\leq i\leq l\}$, $\psi=2x_1=2\sum_{i=1}^{l-1}\g_i+\g_l$ and therefore
$(x_i,x_j)=1/4(\psi,\psi)\de_{ij}$;
by (3.1), we obtain
\begin{eqnarray}
&&e_j=\f{2}{(\psi,\psi)}\sum_{k=1}^j x_k\qquad 1\leq j\leq l;\\
&&(e_i,e_j)=\f{i}{(\psi,\psi)}\qquad 1\leq i\leq j\leq l.
\end{eqnarray}

$\Si=\sw{d}_l$($l\geq 4$): The corresponding Dynkin diagram is
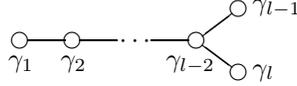
\setlength{\unitlength}{1mm}
\begin{center}
\begin{picture}(33,11.5)
\put(1.5,6){\circle{2}}
\put(2.5,6){\line(1,0){5}}
\put(8.5,6){\circle{2}}
\put(9.5,6){\line(1,0){5}}
\put(14.8,5.1){$\cdots$}
\put(19,6){\line(1,0){5}}
\put(25,6){\circle{2}}
\put(0,2.5){$\g_1$}
\put(7,2.5){$\g_2$}
\put(21,2.5){$\g_{l-2}$}
\put(25.8,6.6){\line(4,3){4}}
\put(30.6,10.2){\circle{2}}
\put(25.8,5.4){\line(4,-3){4}}
\put(30.6,1.8){\circle{2}}
\put(32.5,10){$\g_{l-1}$}
\put(32.5,1.3){$\g_l$}
\end{picture}
\end{center}
Denote $\g_i=x_i-x_{i+1}$($1\leq i\leq l-1$), $\g_l=x_{l-1}+x_l$, then $\Si=\{\pm(x_i\pm x_j):1\leq i\leq j\leq l\}$,
$\psi=x_1+x_2=\g_1+2\sum_{i=2}^{l-2}\g_i+\g_{l-1}+\g_l$ and therefore
$(x_i,x_j)=1/2(\psi,\psi)\de_{ij}$;
by (3.1), we obtain
\begin{eqnarray}
&&e_1=\f{2}{(\psi,\psi)}x_1,\ e_j=\f{1}{(\psi,\psi)}\sum_{k=1}^j x_k\ (2\leq j\leq l-2),\nonumber\\
&&e_{l-1}=\f{1}{(\psi,\psi)}(\sum_{k=1}^{l-1}x_k-x_l),\ e_l=\f{1}{(\psi,\psi)}\sum_{k=1}^l x_k;
\end{eqnarray}
\begin{eqnarray}
&&(e_1,e_1)=\f{2}{(\psi,\psi)},\ (e_1,e_j)=\f{1}{(\psi,\psi)}\ (2\leq j\leq l),\ (e_i,e_j)=\f{i}{2(\psi,\psi)}\ (2\leq i\leq l-2\mbox{ and }j\geq i),\nonumber\\
&&(e_{l-1},e_{l-1})=(e_l,e_l)=\f{l}{2(\psi,\psi)},\ (e_{l-1},e_l)=\f{l-2}{2(\psi,\psi)}.
\end{eqnarray}

$\Si=\sw{e}_6$: The corresponding Dynkin diagram is
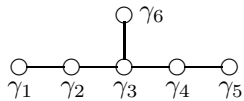
\setlength{\unitlength}{1mm}
\begin{center}
\begin{picture}(31,11.5)
\put(1.5,3.5){\circle{2}}
\put(2.5,3.5){\line(1,0){5}}
\put(8.5,3.5){\circle{2}}
\put(9.5,3.5){\line(1,0){5}}
\put(15.5,3.5){\circle{2}}
\put(16.5,3.5){\line(1,0){5}}
\put(22.5,3.5){\circle{2}}
\put(23.5,3.5){\line(1,0){5}}
\put(29.5,3.5){\circle{2}}
\put(0,0){$\g_1$}
\put(7,0){$\g_2$}
\put(14,0){$\g_3$}
\put(21,0){$\g_4$}
\put(28,0){$\g_5$}
\put(15.5,4.5){\line(0,1){5}}
\put(15.5,10.5){\circle{2}}
\put(17.5,10){$\g_6$}
\end{picture}
\end{center}
Then $\psi=\g_1+2\g_2+3\g_3+2\g_4+\g_5+2\g_6$; since all the roots have the same length,
\begin{eqnarray*}
\Om=\f{1}{2}(\psi,\psi)\left(\begin{array}{cccccc}2&-1& & & &\\-1& 2&-1& & & \\& -1& 2&-1 & & -1\\& & -1& 2& -1& \\& & & -1 & 2& \\ & & -1 & & & 2\end{array}\right);
\end{eqnarray*}
then by Lemma 5.1,
\begin{eqnarray}
\big((e_i,e_j)\big)=\big(\f{1}{d_jd_i}(\Om^{-1})_{ji}\big)=\f{1}{3(\psi,\psi)}\left(\begin{array}{cccccc}8 & 5 & 4 & 4 & 4 & 3\\
                                                                                      5 & 5 & 4 & 4 & 4 & 3\\
                                                                                      4 & 4 & 4 & 4 & 4 & 3\\
                                                                                      4 & 4 & 4 & 5 & 5 & 3\\
                                                                                      4 & 4 & 4 & 5 & 8 & 3\\
                                                                                      3        &  3       &  3       & 3        & 3        & 3
                                                                  \end{array}\right).
\end{eqnarray}

$\Si=\sw{e}_7$: The corresponding Dynkin diagram is
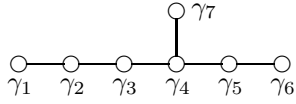
\setlength{\unitlength}{1mm}
\begin{center}
\begin{picture}(38,11.5)
\put(1.5,3.5){\circle{2}}
\put(2.5,3.5){\line(1,0){5}}
\put(8.5,3.5){\circle{2}}
\put(9.5,3.5){\line(1,0){5}}
\put(15.5,3.5){\circle{2}}
\put(16.5,3.5){\line(1,0){5}}
\put(22.5,3.5){\circle{2}}
\put(23.5,3.5){\line(1,0){5}}
\put(29.5,3.5){\circle{2}}
\put(30.5,3.5){\line(1,0){5}}
\put(36.5,3.5){\circle{2}}
\put(0,0){$\g_1$}
\put(7,0){$\g_2$}
\put(14,0){$\g_3$}
\put(21,0){$\g_4$}
\put(28,0){$\g_5$}
\put(35,0){$\g_6$}
\put(22.5,4.5){\line(0,1){5}}
\put(22.5,10.5){\circle{2}}
\put(24.5,10){$\g_7$}
\end{picture}
\end{center}
Then $\psi=\g_1+2\g_2+3\g_3+4\g_4+3\g_5+2\g_6+2\g_7$;
since all the roots have the same length,
\begin{eqnarray*}
\Om=\f{1}{2}(\psi,\psi)\left(\begin{array}{ccccccc}2&-1& & & & &\\-1& 2&-1& & &  &\\& -1 & 2 & -1 & & &\\ & & -1& 2&-1 & & -1\\ & & & -1& 2& -1& \\& & & & -1 & 2& \\& & & -1 & & & 2\end{array}\right);
\end{eqnarray*}
then by Lemma 5.1,
\begin{eqnarray}
\big((e_i,e_j)\big)=\big(\f{1}{d_jd_i}(\Om^{-1})_{ji}\big)=\f{1}{12(\psi,\psi)}\left(\begin{array}{ccccccc}36 & 24 & 20 & 18 & 16 & 12 & 18\\
                                                                                      24 & 24 & 20 & 18 & 16 & 12 & 18\\
                                                                                      20 & 20 & 20 & 18 & 16 & 12 & 18\\
                                                                                      18 & 18 & 18 & 18 & 16 & 12 & 18\\
                                                                                      16 & 16 & 16 & 16 & 16 & 12 & 16\\
                                                                                      12 & 12 & 12 & 12 & 12 & 12 & 12\\
                                                                                      18 & 18 & 18 & 18 & 16 & 12 & 21
                                                                  \end{array}\right).
\end{eqnarray}

From (5.2), (5.8), using the technology given in Section 3, we can give the group structure of $Z_{\wtd{K}}(\wtd{M})$,
which only depend on the type of $\Si$.

\begin{Pro}

(i) If $\Si=\sw{e}_8,\sw{f}_4,\sw{g}_2$ or $(\sw{bc})_l$($l\geq 1$), then $Z_{\wtd{M}}(\wtd{K})=\{\wtd{o}\}=\Z_1$.

(ii) If $\Si=\sw{a}_l$,
then $Z_{\wtd{M}}(\wtd{K})=\big\{\wtd{\Exp}(\pi\sqrt{-1}e_j):1\leq j\leq l)\big\}\cup \{\wtd{o}\}=\Z_{l+1}$,
and $\big(\wtd{\Exp}(\pi\sqrt{-1}e_1)\big)^j=\wtd{\Exp}
(\pi\sqrt{-1}e_j)$ for every $1\leq j\leq l$.

(iii) If $\Si=\sw{b}_l$,
then $Z_{\wtd{M}}(\wtd{K})=\big\{\wtd{\Exp}(\pi\sqrt{-1}e_1),\wtd{o}\big\}=\Z_2$.

(iv) If $\Si=\sw{c}_l$,
then $Z_{\wtd{M}}(\wtd{K})=\big\{\wtd{\Exp}(\pi\sqrt{-1}e_l),\wtd{o}\big\}=\Z_2$.

(v) If $\Si=\sw{d}_l$,
then $Z_{\wtd{M}}(\wtd{K})=\big\{\wtd{\Exp}(\pi\sqrt{-1}e_1),\wtd{\Exp}(\pi\sqrt{-1}e_{l-1}),\wtd{\Exp}(\pi\sqrt{-1}e_l),
\wtd{o}\big\}$. When $l$ is even, it is isomorphic to $\Z_2\oplus \Z_2$; when $l$ is odd, it is isomorphic to
$\Z_4$ and $\big(\wtd{\Exp}(\pi\sqrt{-1}e_{l-1})\big)^2=\wtd{\Exp}(\pi\sqrt{-1}e_{1})$,
$\big(\wtd{\Exp}(\pi\sqrt{-1}e_{l-1})\big)^3=\wtd{\Exp}(\pi\sqrt{-1}e_{l})$.

(vi) If $\Si=\sw{e}_6$,
then $Z_{\wtd{M}}(\wtd{K})=\big\{\wtd{\Exp}(\pi\sqrt{-1}e_1),\wtd{\Exp}(\pi\sqrt{-1}e_5),\wtd{o}\big\}=\Z_3$.

(vii) If $\Si=\sw{e}_7$,
then $Z_{\wtd{M}}(\wtd{K})=\big\{\wtd{\Exp}(\pi\sqrt{-1}e_1),\wtd{o}\big\}=\Z_2$.

\end{Pro}

Proof.
By Proposition 3.1 and (3.7), from the fact
that a group of prime order is a cyclic group, (i),(iii)-(iv), (vi)-(vii) is easily seen.

When $\Si=\sw{a}_l$, by (5.2),
\begin{eqnarray}
&&s_{x_1-x_{j+1}}(e_{j+1}-e_j)=e_1\qquad 1\leq j\leq l-1;\nonumber\\
&&s_{x_1-x_{l+1}}(-e_{l})=e_1.
\end{eqnarray}
Where $s_\g$($\g\in \Si$) the reflection with respect to $\g=0$, which belongs to the Weyl group. (5.12) yields
$\wtd{\Exp}(\pi\sqrt{-1}e_{j+1})\big(\wtd{\Exp}(\pi\sqrt{-1}e_j)\big)^{-1}=\wtd{\Exp}(\pi\sqrt{-1}e_1)$
and $\big(\wtd{\Exp}(\pi\sqrt{-1}e_l)\big)^{-1}=\wtd{\Exp}(\pi\sqrt{-1}e_1)$ and furthermore we have (ii).

When $\Si=\sw{d}_l$, by (5.8),
\begin{eqnarray}
&&s_{x_1+x_2}s_{x_1-x_2}(-e_1)=e_1;\nonumber\\
&&s_{x_2+x_3}s_{x_4+x_5}\cdots s_{x_{l-2}+x_{l-1}}(e_1-e_{l-1})=e_{l}\qquad \mbox{ if }l \mbox{ is even},\nonumber\\
&&s_{x_{l-1}-x_l}s_{x_2+x_3}s_{x_4+x_5}\cdots s_{x_{l-3}+x_{l-2}}(e_1-e_{l-1})=e_{l-1}\qquad \mbox{ if }l\mbox{ is odd}.
\end{eqnarray}
Which implies
\begin{eqnarray}
&&\big(\wtd{\Exp}(\pi\sqrt{-1}e_1)\big)^{-1}=\wtd{\Exp}(\pi\sqrt{-1}e_1);\nonumber\\
&&\wtd{\Exp}(\pi\sqrt{-1}e_1)\big(\wtd{\Exp}(\pi\sqrt{-1}e_{l-1})\big)^{-1}=\left\{
\begin{array}{ll}\wtd{\Exp}(\pi\sqrt{-1}e_l)& l\mbox{ is even};\\
\wtd{\Exp}(\pi\sqrt{-1}e_{l-1}) & l\mbox{ is odd}.\end{array}\right.
\end{eqnarray}
Since $\big|Z_{\wtd{M}}(\wtd{K})\big|=4$, $Z_{\wtd{M}}(\wtd{K})$ is isomorphic to $\Z_2\oplus \Z_2$ or $\Z_4$;
then (v) is easily obtained.
$\Box$

\section{The computation of $i(P_\G)$ and $d(P_\G)$}
\setcounter{equation}{0}

Our assumption and denotation keep invariant. At the beginning of the section, we define two new quantities.

DENOTATION 6.1.
Define
\begin{eqnarray}
i(P_\G)=\min_{x\in P'_\G}(x,x)^{1/2},\qquad d(P_\G)=\max_{x\in P_\G}(x,x)^{1/2}=\max_{x\in P'_\G}(x,x)^{1/2};
\end{eqnarray}
where $(,)$ is an inner product on $\sw{h}_{\sw{p}_0}$ induced by the Killing form on $\sw{g}$.
\newline

In the following we shall compute $i(P_\G)$ and $d(P_\G)$.

By the definition of $P'_\G$, for every $x\in P'_\G$, $(x,\psi)=1$ or $(x,e_j)=1/2(e_j,e_j)$ for some
$j$ such that $\wtd{\Exp}(\pi\sqrt{-1}e_j)\in \G$, which implies $d_j=1$. If $(x,\psi)=1$, then $1=(x,\psi)
\leq (x,x)^{1/2}(\psi,\psi)^{1/2}$, which yields $(x,x)^{1/2}\geq (\psi,\psi)^{-1/2}$; if $(x,e_j)=1/2(e_j,e_j)$,
then $1/2(e_j,e_j)=(x,e_j)\leq (x,x)^{1/2}(e_j,e_j)^{1/2}$, which implies $(x,x)^{1/2}\geq 1/2(e_j,e_j)^{1/2}$.
Thus
\begin{eqnarray}
i(P_\G)\geq \min\Big\{(\psi,\psi)^{-1/2},\f{1}{2}(e_j,e_j)^{1/2}:\wtd{\Exp}(\pi\sqrt{-1}e_j)\in \G\Big\}.
\end{eqnarray}
If the right side of (6.2) is equal to $1/2(e_k,e_k)^{1/2}$ for some $k$, let $x=1/2e_k$, then $(x,\g_i)=1/2\de_{ik}\geq 0$
for every $1\leq i\leq l$, $(x,\psi)=1/2\leq 1$, $(x,e_j)=1/2(e_k,e_j)\leq 1/2(e_k,e_k)^{1/2}(e_j,e_j)^{1/2}\leq 1/2(e_j,e_j)$ for every $j$ such that
$\wtd{\Exp}(\pi\sqrt{-1}e_j)\in \G$; which yields $x\in P'_\G$ and hence $i(P_\G)=1/2(e_k,e_k)^{1/2}$. Otherwise,
the right side of (6.2) is equal to $(\psi,\psi)^{-1/2}$, let $x=\psi/(\psi,\psi)$, then $(x,\g_i)\geq 0$,
$(x,\psi)=1$, $(x,e_j)=(\psi,\psi)^{-1}\leq 1/4(e_j,e_j)\leq 1/2(e_j,e_j)$ for every $j$ such that
$\wtd{\Exp}(\pi\sqrt{-1}e_j)\in \G$; which yields $x\in P'_\G$ and then $i(P_\G)=(\psi,\psi)^{-1/2}$. Therefore
\begin{eqnarray}
i(P_\G)= \min\Big\{(\psi,\psi)^{-1/2},\f{1}{2}(e_j,e_j)^{1/2}:\wtd{\Exp}(\pi\sqrt{-1}e_j)\in \G\Big\}.
\end{eqnarray}

By (6.3) and the results of $(e_i,e_j)$ in Section 5,  we can compute $i(P_\G)$ for any given $\Si$ and $\G$.
We list the results as follows.
\begin{eqnarray}
&&\Si=\sw{a}_l: i(P_\G)=\left\{\begin{array}{ll}\f{\sqrt{2}}{2}(\psi,\psi)^{-1/2}l^{1/2}(l+1)^{-1/2} & \G=\Z_{l+1};\\
(\psi,\psi)^{-1/2}(l-1)^{1/2}(l+1)^{-1/2} & \G=\Z_{\f{l+1}{2}},l\geq 3, l\mbox{ is odd};\\
\f{\sqrt{3}}{2}(\psi,\psi)^{-1/2} & l=5,\G=\Z_2;\\
(\psi,\psi)^{-1/2} & \mbox{otherwise.}
\end{array}\right.\\
&&\Si=\sw{b}_l:i(P_\G)=\f{\sqrt{2}}{2}(\psi,\psi)^{-1/2}\qquad (\G=\Z_2).\\
&&\Si=\sw{c}_l:i(P_\G)=\min\big\{1,\f{1}{2}l^{1/2}\big\}(\psi,\psi)^{-1/2}\qquad (\G=\Z_2).\\
&&\Si=\sw{d}_l:i(P_\G)=\left\{\begin{array}{ll}\f{\sqrt{2}}{2}(\psi,\psi)^{-1/2} & \wtd{\Exp}(\pi\sqrt{-1}e_1)\in \G;\\ \min\big\{1,\f{\sqrt{2}}{4}l^{1/2}\big\}(\psi,\psi)^{-1/2} & \wtd{\Exp}(\pi\sqrt{-1}e_1)\notin \G.\end{array}\right. \\
&&\Si=\sw{e}_6:i(P_\G)=\f{\sqrt{6}}{3}(\psi,\psi)^{-1/2}\qquad (\G=\Z_3).\\
&&\Si=\sw{e}_7:i(P_\G)=\f{\sqrt{3}}{2}(\psi,\psi)^{-1/2}\qquad (\G=\Z_2).
\end{eqnarray}

What about $d(P_\G)$? Notice that $P_\G$ is a convex polyhedron, and for any $x_1,x_2\in P_\G$ and $t\in [0,1]$,
\begin{eqnarray}
\big(tx_1+(1-t)x_2,tx_1+(1-t)x_2\big)^{1/2}&=&\big(t^2(x_1,x_1)+(1-t)^2(x_2,x_2)+2t(1-t)(x_1,x_2)\big)^{1/2}\nonumber\\
&\leq& t(x_1,x_1)^{1/2}+(1-t)(x_2,x_2)^{1/2};
\end{eqnarray}
which yields that $x\in P_\G\mapsto (x,x)^{1/2}$ takes its maximum at the vertices of $P_\G$. It is an elementary
idea to determine all the vertices of $P_\G$ explicitly and then compute $d(P_\G)$. Sometimes the method takes effect,
but when the vertices are too many it doesn't; so we need peculiar tricks for concrete examples. Now we give two lemmas
which will play an important role later.

\begin{Lem} $m\in \Z^+$, $a\geq 0$, $b,s>0$ satisfying $ma\leq s\leq mb$, then if $\la_1,\cdots,\la_m\in [a,b],
\sum_{i=1}^m \la_i\leq s$,
we have
\begin{eqnarray}
\sum_{i=1}^m \la_i^2\leq [\f{mb-s}{b-a}]^2a^2+[\f{s-ma}{b-a}]^2 b^2+c^2; \ \mbox{where }c=s-[\f{mb-s}{b-a}]a-[\f{s-ma}{b-a}]b.
\end{eqnarray}

\end{Lem}

Proof. Denote $D=\{(\la_1,\cdots,\la_m)\in [a,b]^m:\sum_{i=1}^m \la_i\leq s\}$, then $D$ is compact and every
continuous function on $D$ can takes its maximum. Denote by $(\mu_1,\cdots,\mu_m)\in D$ such that
$\sum_{i=1}^m \mu_i^2\geq \sum_{i=1}^m \la_i^2$ for every $(\la_1,\cdots,\la_m)\in D$.
We claim $\{\mu_i:1\leq i\leq m\}\cap (a,b)$ has at most one element. If not, we assume $a<\mu_1\leq \mu_2<b$
without loss of generality, then there exists sufficiently small $\varepsilon>0$ such that
$\mu_1-\varepsilon,\mu_2+\varepsilon\in (a,b)$; let $\la_1=\mu_1-\varepsilon$, $\la_2=\mu_2+\varepsilon$,
$\la_i=\mu_i$ ($3\leq i\leq m$), then obviously $\sum_{i=1}^m \la_i^2>\sum_{i=1}^m \mu_i^2$; which causes a contradiction.
Using the same trick, we can prove $\sum_{i=1}^m \mu_i=s$. Then we can easily obtain (6.11).
$\Box$

\begin{Lem}
Let $m\in \Z^+$, $t_1\geq t_2\geq\cdots\geq t_m\geq 0$, if $\la_1\geq \la_2\geq \cdots\geq \la_m\geq 0$
such that $\sum_{k=1}^j \la_k\leq \sum_{k=1}^j t_k$ for every $1\leq j\leq m$, then $\sum_{j=1}^m \la_j^2\leq \sum_{j=1}^m t_j^2$.
\end{Lem}

Proof. Denote $s_j=\sum_{k=1}^j t_k$, $\mu_j=\sum_{k=1}^j \la_k$ ($1\leq j\leq m$), then $\la_1\geq \la_2\geq \cdots\geq \la_m\geq 0$
and $\sum_{k=1}^j \la_k\leq \sum_{k=1}^j t_k$ if and only if
\begin{eqnarray*}
&&2\mu_1-\mu_2\geq 0,\ 2\mu_2-\mu_1-\mu_3\geq 0,\ \cdots, 2\mu_{m-1}-\mu_{m-2}-\mu_m\geq 0,\mu_m-\mu_{m-1}\geq 0;\\
&&\mu_j\leq s_j\ (1\leq j\leq m).
\end{eqnarray*}
Denote by
$D=\{(z_1,\cdots,z_m)\in \R^m:2z_1-z_2\geq 0, z_m-z_{m-1}\geq 0, 2z_k-z_{k-1}-z_{k+1}\geq 0, z_j\leq s_j \mbox{ for any }
2\leq k\leq m-1\mbox{ and }  1\leq j\leq m\},$
then obviously $D$ is convex and $(\mu_1,\cdots,\mu_m),(s_1,\cdots,s_m)\in D$. Denote
by $\z(t)=(1-t)(\mu_1,\cdots,\mu_m)+t(s_1,\cdots,s_m)$, then $\z(t)\in D$ for $t\in [0,1]$, $\z(0)=(\mu_1,\cdots,\mu_m)$,
$\z(1)=(s_1,\cdots,s_m)$ and $\d{\z}(t)=(s_1-\mu_1,\cdots,s_m-\mu_m)$.

Define $f:D\ra \R$
$$f(z_1,\cdots,z_m)=z_1^2+\sum_{j=1}^{m-1}(z_{j+1}-z_j)^2,$$
then $f(\mu_1,\cdots,\mu_m)=\sum_{j=1}^m \la_j^2$, $f(s_1,\cdots,s_m)=\sum_{j=1}^m t_j^2$; since
$$\f{\p f}{\p z_1}=2(2z_1-z_2)\geq 0,\ \f{\p f}{\p z_j}=2(2z_j-z_{j-1}-z_{j+1})\geq 0(2\leq j\leq m-1),\ \f{\p f}{\p z_m}=2(z_m-z_{m-1})\geq 0,$$
$(f\circ \z)'(t)=\sum_{1\leq j\leq m}(s_j-\mu_j)\f{\p f}{\p z_j}\big(\z(t)\big)\geq 0$ and therefore $f\circ \z(1)\geq f\circ \z(0)$; i.e.,
$\sum_{j=1}^m \la_j^2\leq \sum_{j=1}^m t_j^2$. $\Box$
\newline

In the following we give the detail of computing $d(P_\G)$ for any given $\Si$ and $\G$.
\newline

CASE I. $\Si=\sw{a}_l$ and $\G=\Z_{l+1}$, i.e., $\G=\{\wtd{o},\wtd{\Exp}(\pi\sqrt{-1}e_i):1\leq i\leq l\}$.

The denotation of $x_i$ is similar to Section 5;
from the definition of $P_\G$, by (5.2)-(5.3), $x=\sum_{i=1}^{l+1}\la_i x_i\in P_\G$
($\sum_{i=1}^{l+1}\la_i=0$) if and only if
\begin{eqnarray}
&&\la_1-\la_2\geq 0,\cdots,\la_l-\la_{l+1}\geq 0,\la_1-\la_{l+1}\leq \f{2}{(\psi,\psi)};\nonumber\\
&&(l+1-j)\sum_{k=1}^j \la_k-j\sum_{k=j+1}^{l+1}\la_k\leq \f{j(l+1-j)}{(\psi,\psi)}\qquad (1\leq j\leq l).\nonumber\\
&&\mbox{ i.e., }\sum_{k=1}^j \la_k\leq \f{j(l+1-j)}{(l+1)(\psi,\psi)},\ \sum_{k=j+1}^{l+1} \la_k\geq -\f{j(l+1-j)}{(l+1)(\psi,\psi)}.
\end{eqnarray}
Let $1\leq m\leq l$ such that $\la_m\geq 0$ but $\la_{m+1}<0$. Then by (6.12), we have
\begin{eqnarray*}
\sum_{k=1}^j \la_k\leq \sum_{k=1}^j t_k\ (1\leq j\leq m),\qquad\mbox{ where }t_k=\left\{\begin{array}{ll}\f{l+2-2k}{(l+1)(\psi,\psi)} & 1\leq k\leq [\f{l+1}{2}];\\0 & k>[\f{l+1}{2}].\end{array}\right.
\end{eqnarray*}
Since $t_1\geq t_2\geq \cdots\geq t_m\geq 0$, by Lemma 6.2, we have
\begin{eqnarray}
\sum_{k=1}^m \la_k^2\leq \sum_{k=1}^m t_k^2\leq \sum_{k=1}^{[\f{l+1}{2}]}\f{(l+2-2k)^2}{(l+1)^2(\psi,\psi)^2}=\f{l(l+2)}{6(l+1)(\psi,\psi)^2}.
\end{eqnarray}
On the other hand, from (6.12) we have
\begin{eqnarray*}
&&-\la_{l+1}\geq -\la_{l}\geq \cdots\geq -\la_{m+1}>0;\\
&&-\la_{l+1}\leq t_{l+1},-\la_{l+1}-\la_{l}\leq t_{l+1}+t_{l},\cdots, \sum_{k=m+1}^{l+1} (-\la_k)\leq \sum_{k=m+1}^{l+1} t_k.\\
\end{eqnarray*}
Where
$$t_k=\left\{\begin{array}{ll}\f{2k-l-2}{(l+1)(\psi,\psi)} & k\geq [\f{l+1}{2}];\\ 0 & k<[\f{l+1}{2}].\end{array}\right.$$
Since $t_{l+1}\geq t_{l}\geq \cdots\geq t_{m+1}\geq 0$, by Lemma 5.4, we have
\begin{eqnarray}
\sum_{k=m+1}^{l+1} \la_k^2\leq \sum_{k=m+1}^{l+1} t_k^2\leq \sum_{k=[\f{l+1}{2}]}^{l+1} \f{(2k-l-2)^2}{(l+1)^2(\psi,\psi)^2}=\f{l(l+2)}{6(l+1)(\psi,\psi)^2}.
\end{eqnarray}
(6.13) and (6.14) yield
\begin{eqnarray}
(x,x)=\sum_{k=1}^{l+1} \la_k^2(x_k,x_k)=\sum_{k=1}^{l+1} \la_k^2\cdot \f{1}{2}(\psi,\psi)\leq \f{l(l+2)}{6(l+1)(\psi,\psi)}
\end{eqnarray}
and the equal sign holds if and only if
$$x=\sum_{k=1}^{l+1}\f{l+2-2k}{(l+1)(\psi,\psi)}x_k\in P_\G;$$
so we have
\begin{eqnarray}
d(P_\G)=\f{\sqrt{6}}{6}(\psi,\psi)^{-1/2}l^{1/2}(l+2)^{1/2}(l+1)^{-1/2}.
\end{eqnarray}

CASE II. $\Si=\sw{a}_l$ ($l\geq 3$ is odd) and $\G=\Z_2$, i.e., $\G=\{\wtd{o},\wtd{\Exp}(\pi\sqrt{-1}e_{\f{l+1}{2}})\}$.

At first, notice that the linear automorphism $\varphi$ of $\sw{h}_{\sw{p}_0}$ satisfying $\varphi(\g_i)=\g_{l+1-i}$
keeps $(,)$ invariant, which also satisfies $\varphi(e_i)=e_{l+1-i}$.

The vertices of $\triangle$ are $0,e_1,\cdots,e_l$; by (5.2)-(5.3),
for every $1\leq i<\f{l+1}{2}$, $(e_i,e_{\f{l+1}{2}})\leq 1/2(e_{\f{l+1}{2}},e_{\f{l+1}{2}})$ if and only if
$i\leq\f{l+1}{4}$,
$(e_i,e_{\f{l+1}{2}})\geq 1/2(e_{\f{l+1}{2}},e_{\f{l+1}{2}})$ if and only if $i\geq \f{l+1}{4}$;
so the vertices of $P_\G$ are
\begin{eqnarray}
&&0;\ e_i,\varphi(e_i)\ (1\leq i\leq \f{l+1}{4});\ \f{l+1}{4j}e_j,\f{l+1}{4j}\varphi(e_j)\ (\f{l+1}{4}<j\leq \f{l+1}{2});\nonumber\\
&&v_{i,j},\varphi(v_{i,j}),\ w_{i,j}, \varphi(w_{i,j})\ (1\leq i<\f{l+1}{4}, \f{l+1}{4}<j\leq \f{l+1}{2}).
\end{eqnarray}
Where
\begin{eqnarray}
&&v_{i,j}=\f{1}{4(j-i)}\big((4j-l-1)e_i+(l+1-4i)e_j\big),\nonumber\\
&&w_{i,j}=\f{1}{4(j-i)}\big((4j-l-1)e_i+(l+1-4i)\varphi(e_j)\big).
\end{eqnarray}
By computing, we obtain
\begin{eqnarray}
d(P_\G)=\left\{\begin{array}{ll}\ \ \ \ \ \ \ \ \ (e_{\f{l+1}{4}},e_{\f{l+1}{4}})^{1/2}=\f{\sqrt{6}}{4}(\psi,\psi)^{-1/2}(l+1)^{1/2} & \f{l+1}{2}\mbox{ is even;}\\ (v_{\f{l-1}{4},\f{l+3}{4}},v_{\f{l-1}{4},\f{l+3}{4}})^{1/2}=\f{\sqrt{2}}{4}(\psi,\psi)^{-1/2}(3l-1)^{1/2} & \f{l+1}{2}\mbox{ is odd.}\end{array}\right.
\end{eqnarray}

CASE III. $\Si=\sw{b}_l$ and $\G=\Z_2$, i.e., $\G=\{\wtd{o},\wtd{\Exp}(\pi\sqrt{-1}e_1)\}$.

The vertices of $\triangle$ are $0,e_1,\cdots,e_l$; by (5.5), $(e_1,e_j)=1/2(e_1,e_1)$ for every $2\leq j\leq l$,
so the vertices of $P_\G$ include
$$0,\f{1}{2}e_1,e_2,\cdots,e_l.$$
Since $(1/2 e_1,1/2 e_1)=1/2(\psi,\psi)^{-1}$, $(e_j,e_j)=j/2(\psi,\psi)^{-1}$($2\leq j\leq l$) (by (5.10)), we have
\begin{eqnarray}
d(P_\G)=(e_l,e_l)^{1/2}=\f{\sqrt{2}}{2}(\psi,\psi)^{-1/2}l^{1/2}.
\end{eqnarray}

CASE IV. $\Si=\sw{c}_l$ and $\G=\Z_2$, i.e., $\G=\{\wtd{o},\wtd{\Exp}(\pi\sqrt{-1}e_l)\}$.

By (5.6)-(5.7), $x=\sum_{i=1}^l \la_i x_i\in P_\G$ if and only if
$$\la_1-\la_2\geq 0,\cdots,\la_{l-1}-\la_l\geq 0,\la_l\geq 0,\la_1\leq \f{2}{(\psi,\psi)},\sum_{i=1}^l \la_i\leq \f{l}{(\psi,\psi)}.$$
By Lemma 6.1, if $l$ is even,
$$(x,x)=\sum_{i=1}^l \la_i^2(x_i,x_i)\leq \f{l}{2}\big(\f{2}{(\psi,\psi)}\big)^2\cdot\f{1}{4}(\psi,\psi)=\f{l}{2(\psi,\psi)};$$
if $l$ is odd,
$$(x,x)=\sum_{i=1}^l \la_i^2(x_i,x_i)\leq \Big(\f{l-1}{2}\big(\f{2}{(\psi,\psi)}\big)^2+\big(\f{1}{(\psi,\psi)}\big)^2\Big)\cdot \f{1}{4}(\psi,\psi)=\f{2l-1}{4(\psi,\psi)}.$$
Then
\begin{eqnarray}
d(P_\G)=\left\{\begin{array}{ll}\f{\sqrt{2}}{2}(\psi,\psi)^{-1/2}l^{1/2} & l\mbox{ is even};\\ \f{1}{2}(\psi,\psi)^{-1/2}(2l-1)^{1/2} & l\mbox{ is odd.}\end{array}\right.
\end{eqnarray}

CASE V. $\Si=\sw{d}_l$ and $\G=Z_{\wtd{M}}(\wtd{K})$, i.e., $\G=\{\wtd{o},\wtd{\Exp}(\pi\sqrt{-1}e_i):i=1,l-1\mbox{ or }l\}$.

By (5.8)-(5.9), $x=\sum_{i=1}^l \la_i x_i\in P_{\G}$ if and only if
\begin{eqnarray*}
&&\la_1-\la_2\geq 0,\cdots,\la_{l-1}-\la_l\geq 0,\la_{l-1}+\la_l\geq 0,\la_1+\la_2\leq \f{2}{(\psi,\psi)};\\
&&\la_1\leq \f{1}{(\psi,\psi)}, \sum_{i=1}^{l-1}\la_i-\la_l\leq \f{l}{2(\psi,\psi)},\sum_{i=1}^l \la_i\leq \f{l}{2(\psi,\psi)}.
\end{eqnarray*}
$\la_{l-1}+\la_l\geq 0$ and $\la_{l-1}-\la_l\geq 0$ yield $\la_{l-1}\geq |\la_l|\geq 0$; $\sum_{i=1}^{l-1}\la_i-\la_l\leq \f{l}{2(\psi,\psi)}$
and $\sum_{i=1}^{l}\la_i\leq \f{l}{2(\psi,\psi)}$ yield $\sum_{i=1}^{l-1}\la_i+|\la_l|\leq \f{l}{2(\psi,\psi)}$, then
by Lemma 6.1,
$$(x,x)=\sum_{i=1}^l \la_i^2(x_i,x_i)=(\sum_{i=1}^{l-1}\la_i^2+|\la_l|^2)\cdot \f{1}{2}(\psi,\psi)\leq \left\{\begin{array}{ll}\f{l}{4(\psi,\psi)} & l\mbox{ is even};\\\f{2l-1}{8(\psi,\psi)} & l\mbox{ is odd.}\end{array}\right.$$
and the equal sign holds if and only if $x=\f{1}{(\psi,\psi)}\sum_{i=1}^{\f{l}{2}}x_i\in P_\G$ when $l$ is even,
$x=\f{1}{(\psi,\psi)}\sum_{i=1}^{\f{l-1}{2}}x_i+\f{1}{2(\psi,\psi)}x_{\f{l+1}{2}}\in P_{\G}$ when $l$ is odd.
Thus
\begin{eqnarray}
d(P_\G)=\left\{\begin{array}{ll}\f{1}{2}(\psi,\psi)^{-1/2}l^{1/2} & l\mbox{ is even};\\ \f{\sqrt{2}}{4}(\psi,\psi)^{-1/2}(2l-1)^{1/2} & l\mbox{ is odd}.\end{array}\right.
\end{eqnarray}

CASE VI. $\Si=\sw{d}_l$ and $\G=\{\wtd{o},\wtd{\Exp}(\pi\sqrt{-1}e_1)\}$.

By (5.9), $(e_1,e_j)=1/2(e_1,e_1)$ for every $2\leq j\leq l$,
so the vertices of $P_\G$ include
$$0,\f{1}{2}e_1,e_2,\cdots,e_l.$$
Since $(1/2 e_1,1/2 e_1)=1/2(\psi,\psi)^{-1}$, $(e_j,e_j)=j/2(\psi,\psi)^{-1}$($2\leq j\leq l-2$) and
$(e_{l-1},e_{l-1})=(e_l,e_l)=l/2(\psi,\psi)^{-1}$ (by (5.10)), we have
\begin{eqnarray}
d(P_\G)=(e_{l-1},e_{l-1})^{1/2}=(e_l,e_l)^{1/2}=\f{\sqrt{2}}{2}(\psi,\psi)^{-1/2}l^{1/2}.
\end{eqnarray}

CASE VII. $\Si=\sw{d}_l$ ($l$ is even) and $\G=\{\wtd{o},\wtd{\Exp}(\pi\sqrt{-1}e_{l-1})\}$ or $\{\wtd{o},\wtd{\Exp}(\pi\sqrt{-1}e_l)\}$.

By (5.9),
when $1\leq i\leq \f{l}{2}$, $(e_i,e_{l-1})\leq 1/2(e_{l-1},e_{l-1})$; when $\f{l}{2}\leq i\leq l$, $(e_i,e_{l-1})\geq 1/2(e_{l-1},e_{l-1})$.
Thus if $\G=\{\wtd{o},\wtd{\Exp}(\pi\sqrt{-1}e_{l-1})\}$, the vertices of $P_\G$ are
\begin{eqnarray*}
&&0,\ e_i(1\leq i\leq \f{l}{2}),\ \f{l}{2j}e_j(\f{l}{2}+1\leq j\leq l-2),\ \f{1}{2}e_{l-1},\ \f{l}{2(l-2)}e_l;\\
&&\f{1}{2(l-2-i)}\big((l-4)e_i+(l-2i)e_l\big),\ \f{1}{2(l-i)}\big(le_i+(l-2i)e_{l-1}\big)\ (2\leq i\leq l-2);\\
&&\f{1}{2(j-i)}\big((2j-l)e_i+(l-2i)e_j\big)\ (2\leq i\leq \f{l}{2}-1,\f{l}{2}+1\leq j\leq l-2);\\
&&\f{1}{2}(e_1+e_l),\ \f{1}{2(l-2)}\big(le_1+(l-4)e_{l-1}\big),\ \f{1}{2(j-2)}\big((2j-l)e_1+(l-4)e_j\big)\ (\f{l}{2}+1\leq j\leq l-2).
\end{eqnarray*}
and
\begin{eqnarray}
d(P_\G)=\max\big\{(e_1,e_1)^{1/2},(e_{\f{l}{2}},e_{\f{l}{2}})^{1/2}\big\}=\left\{\begin{array}{cc}
\sqrt{2}(\psi,\psi)^{-1/2} & l\leq 6;\\
\f{1}{2}(\psi,\psi)^{-1/2}l^{1/2} & l\geq 8.
\end{array}\right.
\end{eqnarray}
Similarly, if $\G=\{\wtd{o},\wtd{\Exp}(\pi\sqrt{-1}e_{l})\}$, (6.24) also holds.

CASE VIII. $\Si=\sw{e}_6$ and $\G=\Z_3$, i.e., $\G=\{\wtd{o},\wtd{\Exp}(\pi\sqrt{-1}e_1),\wtd{\Exp}(\pi\sqrt{-1}e_5)\}$.

By (5.10), the vertices of $P_\G$ are
\begin{eqnarray*}
&&0,e_3,e_6,\f{1}{2}e_1,\f{4}{5}e_2,\f{1}{5}e_1+\f{4}{5}e_6,\f{1}{2}e_2+\f{1}{2}e_6,\\
&&\f{4}{5}e_4,\f{1}{2}e_4+\f{1}{2}e_6,\f{2}{3}e_4+\f{1}{6}e_1,\f{4}{9}e_4+\f{4}{9}e_2,\f{4}{9}e_4+\f{1}{9}e_1+\f{4}{9}e_6,\f{1}{3}e_4+\f{1}{3}e_2+\f{1}{3}e_6,\\
&&\f{1}{2}e_5,\f{1}{5}e_5+\f{4}{5}e_6,\f{1}{3}e_5+\f{1}{3}e_1,\f{1}{6}e_5+\f{2}{3}e_2,\f{1}{6}e_5+\f{1}{6}e_1+\f{2}{3}e_6,\f{1}{9}e_5+\f{4}{9}e_2+\f{4}{9}e_6
\end{eqnarray*}
and
\begin{eqnarray}
d(P_\G)=(e_3,e_3)^{1/2}=\f{2\sqrt{3}}{3}(\psi,\psi)^{-1/2}.
\end{eqnarray}

CASE IX. $\Si=\sw{e}_7$ and $\G=\Z_2$, i.e., $\G=\{\wtd{o},\wtd{\Exp}(\pi\sqrt{-1}e_1)\}$.

By (5.11),  the vertices of $P_\G$ are
\begin{eqnarray*}
&&0,e_4,e_5,e_6,e_7,\f{1}{2}e_1,\f{3}{4}e_2,\f{9}{10}e_3,\f{1}{10}e_1+\f{9}{10}e_5,\f{1}{4}e_2+\f{3}{4}e_5,\\
&&\f{1}{2}e_3+\f{1}{2}e_5,\f{1}{4}e_1+\f{3}{4}e_6,\f{1}{2}e_2+\f{1}{2}e_6,\f{3}{4}e_3+\f{1}{4}e_6
\end{eqnarray*}
and
\begin{eqnarray}
d(P_\G)=(e_7,e_7)^{1/2}=\f{\sqrt{7}}{2}(\psi,\psi)^{-1/2}.
\end{eqnarray}

REMARK 6.1.
If $\Si=\sw{a}_l$, $\G=\Z_r$ such that $2<r<l+1$, the author temporarily has no idea to overcome
the difficulty of computing $d(P_\G)$.

\section{The squared length of the highest restricted root}
\setcounter{equation}{0}

Results of this section about $(\psi,\psi)$  compensate Section 6; after computing $(\psi,\psi)$, we can
obtain $i(P_\G)$ and $d(P_\G)$ explicitly.

In this section, we assume $(\sw{u},\th,\lan,\ran)$ be irreducible; the denotation of
$(,),\De,\Si,\sw{g},\sw{h},\sw{h}_\R,\sw{h}_{\sw{p}_0},\De_0,\\m_\g(\g\in \Si)$ is same as Section 2; and denote by $n$ and
$l$ respectively the rank of $\De$ and $\Si$. Then $(\sw{u},\th,\lan,\ran)$ belongs to one of the
two following types: (I) $\sw{u}$ is compact and simple, $\th$ is an involution;
(II) $\sw{u}$ is a product of two compact simple algebras exchanged by $\th$ (See \cite{AB} p.28).
\newline

TYPE I. In the case, $\De$ and $\Si$ are both irreducible; denote by $\de$ the highest root of $\De$;
since the orderings of $\De$ and $\Si$ are compatible, (i.e., $\a\geq \be$ yields $\bar{\a}\geq \bar{\be}$ for
arbitrary $\a,\be\in \De$), $\bar{\de}$ is the highest root of $\Si$, i.e., $\psi=\bar{\de}$.

Denote by $\de^\perp=\{x\in \sw{h}_{\R}:(x,\de)=0\}$, then $\De\cap \de^\perp$ is obviously a subsystem of $\De$
with an induced ordering; let $B=\{\a_1,\cdots,\a_n\}$ be the set of simple roots in $\De$, then
$B\cap \de^\perp$ is the simple root system of $\De\cap \de^\perp$, and $\a_i\in B\cap \de^\perp$ if and only if $\de-\a_i\notin \De\cup \{0\}$;
then according to the Dynkin diagram of $\De$, we can clarify $B\cap \de^\perp$ and $\De\cap \de^\perp$ (for details see \cite{Y}).

On $\De\cap \de^\perp$, we have the following lemmas:

\begin{Lem}
$(\de,\de)=4(|\De|-|\De\cap \de^\perp|+6)^{-1}$.
\end{Lem}

\begin{Lem}
$(\bar{\de},\bar{\de})=(\de,\de)$ or $1/2(\de,\de)$, and the following conditions are
equivalent:

(a) $(\bar{\de},\bar{\de})=(\de,\de)$;

(b) $\de^\th=-\de$;

(c) $B_0\subset B\cap \de^\perp$, where $B_0=B\cap \De_0$;

(d) $m_{\bar{\de}}=1$.
\end{Lem}
For details of the proof of the two Lemmas, see \cite{Y}.

According to Lemma 7.1, from those well known facts of $|\De|$ for every irreducible and reduced root system (see \cite{Hel} pp. 461-474),
we can obtain $(\de,\de)$ as follows:
\begin{eqnarray}
&&\De=\sw{a}_n:(\de,\de)=\f{1}{n+1};\ \De=\sw{b}_n:(\de,\de)=\f{1}{2n-1}; \ \De=\sw{c}_n:(\de,\de)=\f{1}{n+1};\nonumber\\
&&\De=\sw{d}_n:(\de,\de)=\f{1}{2n-2};
\ \De=\sw{e}_6:(\de,\de)=\f{1}{12};\ \De=\sw{e}_7:(\de,\de)=\f{1}{18};\nonumber\\
&&\De=\sw{e}_8:(\de,\de)=\f{1}{30};
\ \De=\sw{f}_4:(\de,\de)=\f{1}{9};\ \De=\sw{g}_2:(\de,\de)=\f{1}{4}.
\end{eqnarray}

By Lemma 7.2, from the Satake diagram given by Araki in \cite{Ar}, we can
justify whether $(\bar{\de},\bar{\de})=(\de,\de)$ or $(\bar{\de},\bar{\de})=1/2(\de,\de)$ for every type of
irreducible, simple and compact orthogonal involutive Lie algebras. The ultimate results are: $(\bar{\de},\bar{\de})=1/2(\de,\de)$ when
$(\sw{u},\th)$ belongs to $A\ II,C\ II,E\ IV,F\ II$ or $(\sw{u},\th)$ belongs to $BD\ I$ and $l=1$;
otherwise $(\bar{\de},\bar{\de})=(\de,\de)$ (for details see \cite{Y}). Combining the results with (7.1), we can compute $(\bar{\de},\bar{\de})$,
i.e., $(\psi,\psi)$.
\newline

TYPE II. In this case, we denote $\sw{u}=\sw{v}\oplus \sw{v}$, where $\sw{v}$ is a compact and simple Lie
algebra; then $\th(X,Y)=(Y,X)$ for arbitrary $X,Y\in \sw{v}$, $\sw{k}_0=\{(X,X):X\in \sw{v}\}$, $\sw{p}_*
=\{(X,-X):X\in \sw{v}\}$. Let $\sw{t}$ be a maximal abelian subalgebra of $\sw{v}$, $\sw{t}_0=\sqrt{-1}\sw{t}$,
$\De^*\subset \sw{t}_0$ be the root system of $\sw{v}\otimes \C$ with respect to $\sw{t}\otimes \C$ with an ordering;
then $\sw{h}_{\sw{p}_*}
=\{(X,-X):X\in \sw{t}\}$ is a maximal abelian space of $\sw{p}_*$ and we can assume $\sw{h}_{\sw{k}_0}=\{(X,X):
X\in \sw{t}\}$; thus $\sw{h}_{\sw{p}_0}=\{(x,-x):x\in \sw{t}_0\}$, $\sw{h}_{\R}=\{(x,y):x,y\in \sw{t}_0\}$ and
\begin{eqnarray}
\De=(\De^*,0)\cup (0,\De^*),\qquad \Si=\big\{(\f{1}{2}\a,-\f{1}{2}\a):\a\in \De^*\big\}.
\end{eqnarray}
$\De$ has an lexicographic ordering induced by the ordering of $\De^*$, and we can define an ordering on $\Si$:
$(1/2\a,-1/2\a)>0$ if and only if $\a>0$; obviously $\De$ and $\Si$ have compatible orderings. Denote by $\de$
the highest root of $\De^*$, then $\psi=(1/2\de,-1/2\de)$ and
\begin{eqnarray}
(\psi,\psi)=\big((\f{1}{2}\de,-\f{1}{2}\de),(\f{1}{2}\de,-\f{1}{2}\de)\big)=\f{1}{2}(\de,\de),
\end{eqnarray}
i.e., the squared length of the highest restricted root is a half of the squared length of the highest root of $\De^*$.

\section{Computation of injectivity radius and diameter}
\setcounter{equation}{0}

From the definition of injectivity radius and diameter of an arbitrary Riemannian manifold, by Theorem 4.1, Denotation 6.1
and (4.4), we have the following Theorem.

\begin{Thm}

Let $(\sw{u},\th,\lan,\ran)$ be a reduced, compact and irreducible orthogonal involutive Lie algebra, $\wtd{M}=\wtd{U}
/\wtd{K}$ be the simply connected Riemannian symmetric space associated with $(\sw{u},\th,\lan,\ran)$,
$M=\wtd{M}/\G$ be a Clifford-Klein form of $\wtd{M}$, where $\G$ is a subgroup of $Z_{\wtd{M}}(\wtd{K})$
satisfying $\G\neq \{\wtd{o}\}$, then
$i(M)=\pi\ep^{1/2}i(P_\G)$ and $d(M)=\pi\ep^{1/2}d(P_\G)$, where $\ep$ is a positive constant such that $\lan,\ran=
-\ep(,)$.

\end{Thm}

REMARK 8.1.
$\ep$ has  geometric meaning. Let $\n$ and $R$ be respectively Levi-Civita connection and curvature tensor on $M$ with respect to the metric $g$
(where $R(X,Y)=-[\n_X,\n_Y]+\n_{[X,Y]}$), then  $R(X,Y)Z=ad[X,Y]Z$ (cf.  \cite{Ko2} p.231, notice the different sign convention
for the curvature tensor); moreover, by choosing an adapted base we have
\begin{eqnarray}
Ric(X,Y)=-\f{1}{2}(X,Y)=\f{1}{2\ep}\lan X,Y\ran
\end{eqnarray}
(cf. \cite{Sa} p.180); i.e., $M$ is an Einstein manifold with Ricci curvature $1/(2\epsilon)$.
\newline

Then from the results obtained in Section 6 and Section 7, we can compute $i(M)$ and $d(M)$ for every type
of non-simply connected, compact and irreducible Riemannian symmetric spaces and list the results in Table 8.1 and Table 8.2.
\newline

\newcommand{\ZZ}[2]{\rule[#1]{0pt}{#2}}
\newcommand{\rb}[1]{\raisebox{#1}[0pt]}
\begin{center}
\begin{tabular}{|c|c|c|c|c|c|}
\multicolumn{6}{c}{\textbf{Table 8.1}}\\
\multicolumn{6}{c}{}\\
\multicolumn{6}{c}{\textit{The injectivity radius and diameter of non-simply connected, compact and irreducible}}\\
\multicolumn{6}{c}{ \textit{Riemannian symmetric spaces of Type I when $\ep=1$, i.e., $Ric=1/2$}}\\
\multicolumn{6}{c}{}\\
\hline
\ZZ{-8pt}{22pt}\textbf{ Type} & $\wtd{M}$ & $\Si$   & $\G$ & $i(M)$  & $d(M)$ \\\hline
                 &                      &                   & \ZZ{-8pt}{22pt}$\Z_n$          & $\f{\sqrt{2}}{2}\pi(n-1)^{1/2}$    &  $\f{\sqrt{6}}{6}\pi(n^2-1)^{1/2}$\\ \cline{4-6}
                 &                      &                   & \rb{-0.5ex}{$\Z_{\f{n}{2}}$}   & \rb{-1.8ex}{$\pi(n-2)^{1/2}$}      &  \rb{-1.8ex}{unknown}\\
                 &                      &                   & $(n\geq 6)$                    &                                    &        \\\cline{4-6}
 $A\ I$          & \rb{1ex}{$SU(n)/SO(n)$}        & $\sw{a}_{n-1}$    &                                & \ZZ{-8pt}{22pt}$\sqrt{2}\pi$($n=4$)               &  $\f{\sqrt{6}}{4}\pi n$ ($4|n$)\\
                 & \rb{1ex}{($n\geq 2$)}          &                   & $\Z_2$                         & $\f{3\sqrt{2}}{2}\pi$($n=6$)                      &  \rb{-0.7ex}{$\f{\sqrt{2}}{4}\pi (3n^2-4n)^{1/2}$}\\
                 &                      &                   &                                & \ZZ{-8pt}{22pt}$\pi n^{1/2}$($n\geq 8$)           &  \rb{0.7ex}{($4\nmid n$)}                       \\\cline{4-6}
                 &                      &                   & \ZZ{-8pt}{22pt}otherwise       & $\pi n^{1/2}$                       &  unknown     \\\hline
                 &                      &                   & \ZZ{-8pt}{22pt}$\Z_n$          & $\sqrt{2}\pi(n-1)^{1/2}$           &  $\f{\sqrt{6}}{3}\pi(n^2-1)^{1/2}$\\ \cline{4-6}
                 &                      &                   & \rb{-0.5ex}{$\Z_{\f{n}{2}}$}   & \rb{-1.8ex}{$2\pi(n-2)^{1/2}$}     &  \rb{-1.8ex}{unknown}\\
                 &                      &                   & $(n\geq 6)$                    &                                    &        \\\cline{4-6}
 $A\ II$         & \rb{1ex}{$SU(2n)/Sp(n)$}       & $\sw{a}_{n-1}$    &                                & \ZZ{-8pt}{22pt}$2\sqrt{2}\pi$($n=4$)              &  $\f{\sqrt{6}}{2}\pi n$ ($4|n$)\\
                 & \rb{1ex}{($n\geq 2$)}          &                   & $\Z_2$                         & $3\sqrt{2}\pi$($n=6$)                             &  \rb{-0.7ex}{$\f{\sqrt{2}}{2}\pi (3n^2-4n)^{1/2}$}\\
                 &                      &                   &                                & \ZZ{-8pt}{20pt}$2\pi n^{1/2}$($n\geq 8$)          &  \rb{0.7ex}{($4\nmid n$)}                       \\\cline{4-6}
                 &                      &                   & \ZZ{-8pt}{22pt}otherwise       & $2\pi n^{1/2}$                                    &  unknown     \\\hline
\end{tabular}
\end{center}

\begin{center}
\begin{tabular}{|c|c|c|c|c|c|}
\multicolumn{6}{c}{\textbf{Table 8.1}(\textit{continued})}\\
\multicolumn{6}{c}{}\\
\hline
\ZZ{-8pt}{22pt}\textbf{ Type} & $\wtd{M}$ & $\Si$   & $\G$ & $i(M)$  & $d(M)$ \\\hline
            & \rb{-2.0ex}{$Gr_{p,p}(\C)$}       &                    &               &           \rb{-2.0ex}{$\f{\sqrt{2}}{2}\pi p$ ($p\leq 3$)}        &     \rb{-0.5ex}{$\pi p$ ($p$ is even)}\\
 $A\ III$   & \rb{-1.0ex}{$(p\geq 2)$}          &  $\sw{c}_p$        &  $\Z_2$       &           \rb{-1.0ex}{$\sqrt{2}\pi p^{1/2}$($p\geq 4$)}          &     $\f{\sqrt{2}}{2}\pi (2p^2-p)^{1/2}$\\
            &                                   &                    &               &                                                                  &     ($p$ is odd)  \\\hline
                      &                              &                                        &                           &    \rb{-1ex}{$\sqrt{3}\pi$ ($n=3$)}       &  \ZZ{-8pt}{20pt}$\f{\sqrt{2}}{2}\pi(n^2+n)^{1/2}$\\
\rb{-1ex}{$C\ I$}     &  $Sp(n)/U(n)$                &      \rb{-1ex}{$\sw{c}_n$}             &     \rb{-1ex}{$\Z_2$}     &                                           &  ($n$ is even)\\
                      &  ($n\geq 3$)                 &                                        &                           &    \rb{1ex}{$\pi(n+1)^{1/2}$}             &  $\f{1}{2}\pi(2n^2+n-1)^{1/2}$\\
                      &                              &                                        &                           &    \rb{1ex}{($n\geq 4$)}                  &  ($n$ is odd)   \\\hline
                      &                              &                                        &                           &    $\f{\sqrt{2}}{2}\pi (2p^2+p)^{1/2}$        &  \ZZ{-8pt}{20pt}$\pi(2p^2+p)^{1/2}$\\
\rb{-1ex}{$C\ II$}    &  {$Gr_{p,p}(\H)$}            &      \rb{-1ex}{$\sw{c}_n$}             &     \rb{-1ex}{$\Z_2$}     &    $(p\leq 3)$                                    &  ($n$ is even)\\
                      &  ($p\geq 2$)                 &                                        &                           &    $\sqrt{2}\pi(2p+1)^{1/2}$                      &  $\f{\sqrt{2}}{2}\pi(4p^2-1)^{1/2}$\\
                      &                              &                                        &                           &    ($p\geq 4$)                                    &  ($n$ is odd)   \\\hline
                 &   \rb{-0.3ex}{$Gr_{p,q}(\R)$}        &  \rb{-1.8ex}{$\sw{b}_p$}     &  \rb{-1.8ex}{$\Z_2$}     &  \rb{-1.8ex}{$\f{\sqrt{2}}{2}\pi(p+q-2)^{1/2}$} &  \rb{-1.8ex}{$\f{\sqrt{2}}{2}\pi(p^2+pq-2p)^{1/2}$}\\
                 &   \rb{1.3ex}{$(2<p<q)$}              &                              &                          &                                                 &                                        \\\cline{2-6}
                 &   \ZZ{-8pt}{22pt}$S^q$               &       $\sw{a}_1$             &    $\Z_2$     &  $\f{\sqrt{2}}{2}\pi(q-1)^{1/2}$   &   $\f{\sqrt{2}}{2}\pi(q-1)^{1/2}$\\\cline{2-6}
                 &                              &                              &                                &                         &   \ZZ{-8pt}{20pt}$\f{\sqrt{2}}{2}\pi(p^2-p)^{1/2}$\\
                 &                              &                              &   $Z_{\wtd{M}}(\wtd{K})$       &   $\pi(p-1)^{1/2}$      &   \rb{1.0ex}{($p$ is even)}\\
$BD\ I$          &                              &                              &                                &                         &   $\f{1}{2}\pi(2p^2-3p+1)^{1/2}$\\
                 &                              &                              &                                &                         &   ($p$ is odd)\\\cline{4-6}
                 & $Gr_{p,p}(\R)$               &    $\sw{d}_p$                         &  \ZZ{-8pt}{22pt}$\{\wtd{\Exp}(\pi\sqrt{-1}e_1),\wtd{o}\}$     &   $\pi(p-1)^{1/2}$                 &  $\pi(p^2-p)^{1/2}$ \\\cline{4-6}
                 & ($p\geq 4$)                  &                                       &  \rb{-1ex}{$\{\wtd{\Exp}(\pi\sqrt{-1}e_{p-1}),\wtd{o}\}$}                &   \rb{-0.5ex}{$\f{1}{2}\pi (p^2-p)^{1/2}$} &  \rb{-0.5ex}{$2\pi(p-1)^{1/2}$}\\
                 &                              &                                       &  \rb{-1ex}{$\{\wtd{\Exp}(\pi\sqrt{-1}e_p),\wtd{o}\}$}                    &   ($p\leq 6$)                      &  ($p\leq 6$)\\
                 &                              &                                       &  \rb{-1ex}{($p$ is even)}                                                &   $\sqrt{2}\pi(p-1)^{1/2}$         &  $\f{\sqrt{2}}{2}\pi(p^2-p)^{1/2}$  \\
                 &                              &                                       &                                                                          &   \rb{0.5ex}{($p\geq 8$)}          &  \rb{0.5ex}{($p\geq 8$)}         \\\hline
                      &                              &                                        &                           &    $\f{1}{2}\pi(n^2-n)^{1/2}$        &  \ZZ{-8pt}{20pt}$\f{\sqrt{2}}{2}\pi(n^2-n)^{1/2}$\\
\rb{-1ex}{$D\ III$}   &  $SO(2n)/U(n)$               &      \rb{-1ex}{$\sw{c}_{\f{n}{2}}$}    &     \rb{-1ex}{$\Z_2$}     &    ($n\leq 6$)                       &  ($4|n$)\\
                      &  ($n\geq 4$ is even)         &                                        &                           &    $\sqrt{2}\pi(n-1)^{1/2}$          &  $\f{\sqrt{2}}{2}\pi(n-1)$\\
                      &                              &                                        &                           &    ($n\geq 8$)                       &  ($4\nmid n$)   \\\hline
\ZZ{-8pt}{22pt}$E\ I$    &  $(\sw{e}_6,sp(4))$   &    $\sw{e}_6$       &     $\Z_3$     &   $2\sqrt{2}\pi$    &   $4\pi$\\ \hline
\ZZ{-8pt}{22pt}$E\ IV$   &  $(\sw{e}_6,\sw{f}_4)$&    $\sw{a}_2$       &     $\Z_3$     &   $2\sqrt{2}\pi$    &   $\f{4\sqrt{6}}{3}\pi$\\ \hline
\ZZ{-8pt}{22pt}$E\ V$    &  $(\sw{e}_7,su(8)$    &    $\sw{e}_7$       &     $\Z_2$     &   $\f{3\sqrt{6}}{2}\pi$  & $\f{3\sqrt{14}}{2}\pi$\\ \hline
\ZZ{-8pt}{22pt}$E\ VII$  &  $(\sw{e}_7,\sw{e}_6\oplus \R)$ & $\sw{c}_3$&     $\Z_2$     &   $\f{3\sqrt{6}}{2}\pi$  & $3\sqrt{3}\pi$\\\hline
\end{tabular}
\end{center}

\begin{center}
\begin{tabular}{|c|c|c|c|c|}
\multicolumn{5}{c}{\textbf{Table 8.2}}\\
\multicolumn{5}{c}{}\\
\multicolumn{5}{c}{\textit{The injectivity radius and diameter of non-simply connected, compact and irreducible}}\\
\multicolumn{5}{c}{\textit{Riemannian symmetric spaces of Type II when $\ep=1$, i.e., $Ric=1/2$}}\\
\multicolumn{5}{c}{}\\
\hline
\ZZ{-8pt}{22pt} $\wtd{M}$ & $\De^*$   & $\G$ & $i(M)$  & $d(M)$ \\\hline
                        &                   & \ZZ{-8pt}{22pt}$\Z_n$          & $\pi(n-1)^{1/2}$                   &  $\f{\sqrt{3}}{3}\pi(n^2-1)^{1/2}$\\ \cline{3-5}
                        &                   & \rb{-0.5ex}{$\Z_{\f{n}{2}}$}   & \rb{-1.8ex}{$\sqrt{2}\pi(n-2)^{1/2}$}     &  \rb{-1.8ex}{unknown}\\
                        &                   & $(n\geq 6)$                    &                                    &        \\\cline{3-5}
\rb{1ex}{$SU(n)$}       & $\sw{a}_{n-1}$    &                                & \ZZ{-8pt}{22pt}$2\pi$($n=4$)               &  $\f{\sqrt{3}}{2}\pi n$ ($4|n$)\\
\rb{1ex}{($n\geq 2$)}   &                   & $\Z_2$                         & $3\pi$($n=6$)                              &  \rb{-0.7ex}{$\f{1}{2}\pi (3n^2-4n)^{1/2}$}\\
                        &                   &                                & \ZZ{-8pt}{22pt}$\sqrt{2}\pi n^{1/2}$($n\geq 8$)           &  \rb{0.7ex}{($4\nmid n$)}                       \\\cline{3-5}
                        &                   & \ZZ{-8pt}{22pt}otherwise       & $\sqrt{2}\pi n^{1/2}$                       &  unknown     \\\hline
\rb{-0.3ex}{$Spin(2n+1)$}        &  \rb{-1.8ex}{$\sw{b}_n$}     &  \rb{-1.8ex}{$\Z_2$}     &  \rb{-1.8ex}{$\pi(2n-1)^{1/2}$} &  \rb{-1.8ex}{$\pi(2n^2-n)^{1/2}$}\\
\rb{1.3ex}{$(n\geq 2)$}          &                              &                          &                                                 &                                        \\\hline
                      &                                        &                           &    \rb{-1ex}{$\sqrt{6}\pi$ ($n=3$)}      &  \ZZ{-8pt}{20pt}$\pi(n^2+n)^{1/2}$\\
$Sp(n)$               &      \rb{-1ex}{$\sw{c}_n$}             &     \rb{-1ex}{$\Z_2$}     &                                           &  ($n$ is even)\\
 ($n\geq 3$)          &                                        &                           &    \rb{1ex}{$\sqrt{2}\pi(n+1)^{1/2}$}            &  $\f{\sqrt{2}}{2}\pi(2n^2+n-1)^{1/2}$\\
                      &                                        &                           &    \rb{1ex}{($n\geq 4$)}                  &  ($n$ is odd)   \\\hline
                 &                              &                                &                         &   \ZZ{-8pt}{20pt}$\pi(n^2-n)^{1/2}$\\
                 &                              &   $Z_{\wtd{M}}(\wtd{K})$       &   $\sqrt{2}\pi(n-1)^{1/2}$     &   \rb{1.0ex}{($n$ is even)}\\
                 &                              &                                &                         &   $\f{\sqrt{2}}{2}\pi(2n^2-3n+1)^{1/2}$\\
                 &                              &                                &                         &   ($n$ is odd)\\\cline{3-5}
$Spin(2n)$       &    $\sw{d}_n$                         &  \ZZ{-8pt}{22pt}$\{\wtd{\Exp}(\pi\sqrt{-1}e_1),\wtd{o}\}$     &   $\sqrt{2}\pi(n-1)^{1/2}$                 &  $\sqrt{2}\pi(n^2-n)^{1/2}$ \\\cline{3-5}
($n\geq 4$)      &                                       &  \rb{-1ex}{$\{\wtd{\Exp}(\pi\sqrt{-1}e_{n-1}),\wtd{o}\}$}                &   \rb{-1ex}{$\f{\sqrt{2}}{2}\pi (n^2-n)^{1/2}$} &  \rb{-1ex}{$2\sqrt{2}\pi(n-1)^{1/2}$}\\
                 &                                       &  \rb{-1ex}{$\{\wtd{\Exp}(\pi\sqrt{-1}e_n),\wtd{o}\}$}                    &   ($n\leq 6$)                                   &  ($n\leq 6$)\\
                 &                                       &  \rb{-1ex}{($n$ is even)}                                                &   $2\pi(n-1)^{1/2}$                             &  $\pi(n^2-n)^{1/2}$  \\
                 &                                       &                                                                          &   \rb{0.5ex}{($n\geq 8$)}                       &  \rb{0.5ex}{($n\geq 8$)} \\\hline
\ZZ{-8pt}{22pt}$E_6$  &  $\sw{e}_6$        &   $\Z_3$     &   $4\pi$     &   $4\sqrt{2}\pi$\\\hline
\ZZ{-8pt}{22pt}$E_7$  &  $\sw{e}_7$        &   $\Z_2$     &   $3\sqrt{3}\pi$     &   $3\sqrt{7}\pi$\\\hline
\end{tabular}
\end{center}

REMARK 8.2.
In Table 8.1, $M=\wtd{M}/\G$, where $\wtd{M}$ is the universal covering space of $M$, $\G$ is a subgroup of
$Z_{\wtd{M}}(\wtd{K})=\{p\in \wtd{M}:\tau(k)p=p\mbox{ for every }k\in \wtd{K}\}$; $\Si$ denotes the restricted
root system; $i(M)$ and $d(M)$ are respectively the injective diameter and the diameter of $M$.
In Table 8.2, $M$ is a non-simply connected, compact and simple Lie group with bi-invariant metric
and $\wtd{M}$ is the universal covering group of $M$ with pullback metric; in this case, $Z_{\wtd{M}}(\wtd{K})$
coincides with the center of $\wtd{M}$; let $\sw{v}$ be the Lie algebra associated to $M$, $\sw{t}$ be a maximal abelian
subalgebra of $\sw{v}$,
then $\De^*$ denotes the root system of $\sw{v}\otimes \C$ with respect to $\sw{t}\otimes \C$ (cf. Section 7).
\newline

REMARK 8.3.
In Table 8.1, we identify $\sw{b}_2$ and $\sw{c}_2$.
\newline

REMARK 8.4.
When $\wtd{M}=Gr_{p,p}(\R)$, $\De=\Si=\sw{d}_p$; the Satake diagram of $(B,\th)$ is
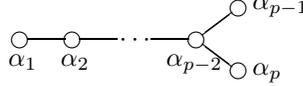
\begin{center}
\begin{picture}(33,11.5)
\put(1.5,6){\circle{2}}
\put(2.5,6){\line(1,0){5}}
\put(8.5,6){\circle{2}}
\put(9.5,6){\line(1,0){5}}
\put(14.8,5.1){$\cdots$}
\put(19,6){\line(1,0){5}}
\put(25,6){\circle{2}}
\put(0,2.5){$\a_1$}
\put(7,2.5){$\a_2$}
\put(21,2.5){$\a_{p-2}$}
\put(25.8,6.6){\line(4,3){4}}
\put(30.6,10.2){\circle{2}}
\put(25.8,5.4){\line(4,-3){4}}
\put(30.6,1.8){\circle{2}}
\put(32.5,10){$\a_{p-1}$}
\put(32.5,1.3){$\a_p$}
\end{picture}
\end{center}
(cf. \cite{Ar}) and the Dynkin
diagram of $\Si$ is
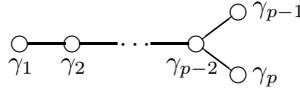
\begin{center}
\begin{picture}(33,11.5)
\put(1.5,6){\circle{2}}
\put(2.5,6){\line(1,0){5}}
\put(8.5,6){\circle{2}}
\put(9.5,6){\line(1,0){5}}
\put(14.8,5.1){$\cdots$}
\put(19,6){\line(1,0){5}}
\put(25,6){\circle{2}}
\put(0,2.5){$\g_1$}
\put(7,2.5){$\g_2$}
\put(21,2.5){$\g_{p-2}$}
\put(25.8,6.6){\line(4,3){4}}
\put(30.6,10.2){\circle{2}}
\put(25.8,5.4){\line(4,-3){4}}
\put(30.6,1.8){\circle{2}}
\put(32.5,10){$\g_{p-1}$}
\put(32.5,1.3){$\g_p$}
\end{picture}
\end{center}
where $\g_i=\bar{\a}_i$ ($1\leq i\leq p$); furthermore, since $\sw{h}_\R=\sw{h}_{\sw{p}_0}$, we have
$\th(\a_i)=-\a_i$ and $\g_i=\a_i$. The definition of $e_1,\cdots, e_{p-1},e_p$ is similar to (3.1).
Let $\phi$ be a linear automorphism of $\sw{h}_\R$ such that $\phi(\a_i)=\a_i$($1\leq i\leq p-2$), $\phi(\a_{p-1})=\a_p$
and $\phi(\a_p)=\a_{p-1}$, then $\phi$ keeps $(,)$ invariant and can be extended to an automorphism of $so(2p)$, which is also denoted by
$\phi$; since $\phi$ commutes with $\th$, it induces an isometry $F$ of $Gr_{p,p}(\R)$, which satisfies
\begin{eqnarray}
&&F\big(\wtd{\Exp}(\pi\sqrt{-1}e_i)\big)=\wtd{\Exp}(\pi\sqrt{-1}e_i)\ (1\leq i\leq p-2);\nonumber\\
&&F\big(\wtd{\Exp}(\pi\sqrt{-1}e_{p-1})\big)=\wtd{\Exp}(\pi\sqrt{-1}e_p),\ F\big(\wtd{\Exp}(\pi\sqrt{-1}e_p)\big)=\wtd{\Exp}(\pi\sqrt{-1}e_{p-1}).
\end{eqnarray}
So when $p$ is even, $Gr_{p,p}(\R)/\{\wtd{\Exp}(\pi\sqrt{-1}e_{p-1}),\wtd{o}\}$ and $Gr_{p,p}(\R)/\{\wtd{\Exp}(\pi\sqrt{-1}e_p),
\wtd{o}\}$ are isometric to each other. Especially, when $p=4$, an arbitrary linear automorphism $\phi$ of $\sw{h}_\R$
satisfying $\phi(B)=B$ and $\phi(\a_2)=\a_2$ keeps $(,)$ invariant, which yields that $Gr_{4,4}(\R)/\{\wtd{\Exp}(\pi\sqrt{-1}e_1),\wtd{o}\}$,
$Gr_{4,4}(\R)/\{\wtd{\Exp}(\pi\sqrt{-1}e_3),\wtd{o}\}$, $Gr_{4,4}(\R)/\{\wtd{\Exp}(\pi\sqrt{-1}e_4),\wtd{o}\}$ are isometric
to each other. On the other hand, if $p$ is even and $p\geq 6$, then $M_1=Gr_{p,p}(\R)/\{\wtd{\Exp}(\pi\sqrt{-1}e_1),\wtd{o}\}$
isn't isometric to $M_2=Gr_{p,p}(\R)/\{\wtd{\Exp}(\pi\sqrt{-1}e_{p-1}),\wtd{o}\}$, although the fundamental
group of them are both isomorphic to $\Z_2$; it is easily seen from Table 8.1 (since $i(M_1)\neq i(M_2)$, $d(M_1)\neq d(M_2)$).

Similarly,
$Spin(8)/\{\wtd{\Exp}(\pi\sqrt{-1}e_1),\wtd{o}\}$ (i.e., $SO(8)$), $Spin(8)/\{\wtd{\Exp}(\pi\sqrt{-1}e_3),\wtd{o}\}$,
$Spin(8)/\\\{\wtd{\Exp}(\pi\sqrt{-1}e_4),\wtd{o}\}$ are isometric to each other; if $n$ is even and $n\geq 6$,
$Spin(2n)/\{\wtd{\Exp}(\pi\sqrt{-1}e_{n-1}),\wtd{o}\}$ and $Spin(2n)/\{\wtd{\Exp}(\pi\sqrt{-1}e_n),\wtd{o}\}$
are isometric to each other, but both of them aren't isometric to $Spin(2n)/\{\wtd{\Exp}(\pi\sqrt{-1}e_1),\wtd{o}\}$
(i.e., $SO(2n)$).
\newline

REMARK 8.5.
In Table 8.1 and Table 8.2, we assume $\ep=1$, i.e., the $K$-invariant metric on $M=U/K$ is induced by $-(,)$ on
$\sw{u}$, and $Ric=1/2$. For general cases such that $\ep\neq 1$, we should multiply the corresponding results in Table 8.1
or Table 8.2 by $\ep^{1/2}$.

For example, let $M=\R\P^q=S^q/\Z_2$ with canonical metric $g$ such that $K=1$, then $Ric=q-1$ and Remark 8.1
yields $\ep=1/\big(2(q-1)\big)$; according to Table 8.1,
\begin{eqnarray}
i(M)=d(M)=\f{\sqrt{2}}{2}\pi(q-1)^{1/2}\ep^{1/2}=\f{\pi}{2}.
\end{eqnarray}
The result is well-known.
\newline

ACKNOWLEDGEMENT. The author wishes to express his sincere gratitude to his supervisor, Professor Y.L. Xin , for his inspiring suggestions, as well as to
Doctor X.S. Liu for providing some references.

\begin{flushright}
Doctor 05 Grade 2\\
Institute of Mathematics\\
Fudan University\\
Shanghai, 200433, P.R. China\\
email: 051018016@fudan.edu.cn\\
\end{flushright}


\begin{thebibliography}{AB}
\bibitem{Ar}Araki, S. I: \textit{On root systems and an infinitesimal classification of irreducible symmetric spaces}, J.Math. Osaka City Univ.\textbf{13}(1962), 1-34.
\bibitem{AB}A. Borel: Semisimple Groups and Riemannian Symmetric Spaces, Hindustan Book Agency, 1998.
\bibitem{CE}J. Cheeger and D.Ebin: Comparison Theorem in Riemannian Geometry, Noth-Holland publishing company, 1975.
\bibitem{Ch}J. Cheeger: \textit{Pinching theorems for a certain class of Riemannian manifolds}, Amer. J. Math. \textbf{91}(1969),807-834.
\bibitem{Cr}R. Crittenden: \textit{Minimum and conjugate points in symmetric spaces}, Canadian J. Math. \textbf{14}(1962),320-328.
\bibitem{GK}M. Goto and E.T. Kobayashi: \textit{On the subgroups of the centers of simiply connected simple Lie groups
-classification of simple Lie groups in the large}, Osaka J. Math. \textbf{6}(1969), 251-281.
\bibitem{GS}K. Grove and K. Shiohama:
\textit{A generalized sphere theorem,}
Ann. Math. (2) \textbf{106} (1977), 201-211.
\bibitem{Hel} S. Helgason: Differential Geometry, Lie groups, and Symmetric Spaces, Academic Press, 1978.
\bibitem{Ic} R. Ichida: \textit{The injectivity radius and the fundamental group of compact homogeneous Riemannian manifolds
of positive curvature,} Tsukuba J. Math. \textbf{24} (2000), 139-156.
\bibitem{Kna} A.W.Knapp: Lie group, Beyond an introduction, Progress in mathematics, 2002.
\bibitem{Ko}S. Kobayashi and K. Nomizu: Foundations of differential geometry, Volume 1, Interscience Publishers, 1963.
\bibitem{Ko2}S. Kobayashi and K. Nomizu: Foundations of differential geometry, Volume 2, Interscience Publishers, 1963.
\bibitem{L}X.S. Liu: \textit{Curvature Estimates for Irreducible Symmetric Spaces}, Chinese Ann. Math. Ser. B \textbf{27} (2006), 287-302.
\bibitem{P} T. P\"{u}ttmann: \textit{Injectivity radius and diameter of the manifolds of flags in the projective planes,}  Math. Z. \textbf{246} (2004), 795-809.
\bibitem{RH}R.Bott and H.Samelson: \textit{Application of the theory of Morse to symmetric spaces}, Amer. J. Math. \textbf{80}(1958), 964-1029.
\bibitem{Sa}Sakai, T.: Riemannian Geometry, Translations of Mathematical Monographs, Vol.149, 1996.
\bibitem{Su}K. Sugahara: \textit{On the cut locus and the topology of Riemannian manifolds,} J. Math. Kyoto Univ.
\textbf{14} (1974), 391-411.
\bibitem{Ta}M. Takeuchi: \textit{On the fundamental group and the group of isometries of a symmetric space,} J. Fac. Sci. Univ. Tokyo Sect. I \textbf{10} (1964), 88-123.
\bibitem{W}A. Weinstein: \textit{The cut locus and conjugate locus of a riemannian manifold}, Ann. of Math. 2nd Ser. \textbf{87}
(1968),29-41.
\bibitem[20]{Y}Ling Yang: \textit{Injectivity radius and Cartan polyhedron for simply connected symmetric spaces,} http://arxiv.org/PS\_cache/math/ps/0609/0609627.ps.gz.

\end{thebibliography}
\end{document}